\documentclass[a4paper,11pt,reqno]{amsart}
\usepackage[english]{babel}
\usepackage{color,amsmath, amssymb, latexsym,mathptmx, amsfonts}

\usepackage[colorlinks=true]{hyperref} 
\usepackage[usenames,dvipsnames]{xcolor}

\usepackage[matrix,arrow,curve,frame]{xy}    

\xymatrixcolsep{1.9pc}                          
\xymatrixrowsep{1.9pc}
\newdir{ >}{{}*!/-5pt/\dir{>}}                  

 \calclayout

\DeclareMathOperator{\nis}{nis}

\DeclareMathOperator{\Hom}{Hom}
\DeclareMathOperator{\colim}{colim}
\DeclareMathOperator{\spec}{Spec}

\renewcommand{\leq}{\leqslant}
\renewcommand{\geq}{\geqslant}
\renewcommand{\phi}{\varphi}

\renewcommand{\kappa}{\varkappa}

\newtheorem{theorem}{Theorem}[section]
\newtheorem*{theorema}{Theorem A}
\newtheorem*{theoremb}{Theorem B}
\newtheorem*{theoremc}{Theorem C}

\newtheorem*{theoremd}{Theorem D}

\newtheorem{lemma}[theorem]{Lemma}
\newtheorem{proposition}[theorem]{Proposition}
\newtheorem{corollary}[theorem]{Corollary}
\newtheorem{definition}[theorem]{Definition}

\newtheorem{remark}[theorem]{Remark}

\newcommand{\bl}[1]{\buildrel #1\over}
\newcommand{\bb}{\mathbb}
\newcommand{\A}{\mathbb{A}}
\newcommand{\Spec}{\operatorname{Spec}}
\newcommand{\pt}{{\rm pt}}
\newcommand{\Gm}{{\mathbb{G}_m}}
\newcommand{\ZF}{\operatorname{\mathbb{Z}F}}
\newcommand{\ZFr}{\operatorname{\mathbb{Z}Fr}}
\newcommand{\Fr}{\operatorname{Fr}}
\newcommand{\id}{\operatorname{id}}

\newcommand{\cc}{\mathcal}
\newcommand{\op}{{\textrm{\rm op}}}

\begin{document}
\sloppy

\title{Cancellation theorem for framed motives of algebraic varieties}
\author{A. Ananyevskiy}
\address{St. Petersburg Branch of Steklov Mathematical Institute, Fontanka
27, St. Petersburg, 191023, Russia}

\email{alseang@gmail.com}

\author{G. Garkusha}
\address{Department of Mathematics, Swansea University, Fabian Way, Swansea SA1 8EN, United Kingdom}
\email{g.garkusha@swansea.ac.uk}

\author{I. Panin}
\address{St. Petersburg Branch of Steklov Mathematical Institute,
Fontanka 27, St. Petersburg, 191023, Russia}
\email{paniniv@gmail.com}




\begin{abstract}
The machinery of framed (pre)sheaves was developed by
Voevodsky~\cite{V2}. Based on the theory, framed motives of
algebraic varieties are introduced and studied in~\cite{GP1}. An
analog of Voevodsky's Cancellation Theorem~\cite{V1} is proved in
this paper for framed motives stating that a natural map of framed
$S^1$-spectra
   $$M_{fr}(X)(n)\to\underline{\Hom}(\mathbb G,M_{fr}(X)(n+1)),\quad n\geq 0,$$
is a schemewise stable equivalence, where $M_{fr}(X)(n)$ is the
$n$th twisted framed motive of $X$. This result is also necessary
for the proof of the main theorem of~\cite{GP1} computing fibrant
resolutions of suspension $\bb P^1$-spectra $\Sigma^\infty_{\bb
P^1}X_+$ with $X$ a smooth algebraic variety.

The Cancellation Theorem for framed motives is reduced to the
Cancellation Theorem for linear framed motives stating that the
natural map  of complexes of abelian groups
\[
\ZF(\Delta^\bullet \times X,Y) \to \ZF((\Delta^\bullet \times
X)\wedge (\Gm,1),Y\wedge (\Gm,1)),\quad X,Y\in Sm/k,
\]
is a quasi-isomorphism, where $\ZF(X,Y)$ is the group of stable
linear framed correspondences in the sense of~\cite{GP1}.
\end{abstract}

\maketitle

\thispagestyle{empty} \pagestyle{plain}

\section{Introduction}

The main goal of the Voevodsky theory on framed correspondences
(see~\cite[Introduction]{V2}) is to suggest a new approach to the
stable motivic homotopy theory $SH(k)$  over a field $k$. This
approach is more amenable to explicit calculations. Recall that
Voevodsky~\cite[Section~2]{V2} invented a category of framed
correspondences $Fr_*(k)$ whose objects are those of $Sm/k$ and
morphisms sets $Fr_*(X,Y)=\sqcup_{n\geq 0} Fr_n(X,Y)$ are defined by
means of certain geometric data. The elements of $Fr_n(X,Y)$ are
called {\it framed correspondences of level $n$}. The definition of
$Fr_*(k)$ is recalled in Section~\ref{preliminaries} below. For
every $Y\in Sm/k$ there is a distinguished morphism
$\sigma_Y=(Y\times 0, Y\times \bb A^1,t,pr_Y) \in Fr_1(Y,Y)$.
Following Voevodsky~\cite{V2}, we denote by
\[
Fr(X,Y):=\colim(Fr_0(X,Y)\xrightarrow{\sigma_Y}Fr_1(X,Y)\xrightarrow{\sigma_Y}\dots
\xrightarrow{\sigma_Y} Fr_n(X,Y)\xrightarrow{\sigma_Y}\dots)
\]
and refer to it as the \textit{set of stable framed
correspondences}. Replacing $Y$ by a simplicial object $Y^\bullet$
in $Sm/k$, we get a simplicial set $Fr(X,Y^\bullet)$. Finally, one
can take the diagonal of the pointed bisimplicial set
$Fr(\Delta^\bullet \times X,Y^\bullet)$. Voevodsky conjectured that
if the motivic space $Fr(\Delta^\bullet \times -,Y^\bullet)$ is
locally connected in the Nisnevich topology, then it is isomorphic
in $H_{\mathbb A^1}(k)$ to the motivic space
$\Omega^{\infty}_{\mathbb P^1}\Sigma^{\infty}_{\mathbb
P^1}(Y^\bullet_+)$. In particular, the theory of framed
correspondences gives a machinery for computing motivic infinite
loop spaces.

Inspired by the Voevodsky theory~\cite{V2}, the theory of framed
motives of algebraic varieties is introduced and developed
in~\cite{GP1}. As an application, the above Voevodsky
conjecture is solved in~\cite[Section~10]{GP1} in the affirmative.
Moreover, under the above assumption on
$Y^\bullet$ the motivic space $Fr(\Delta^\bullet \times
-,Y^\bullet)$ is $\bb A^1$-local. This result can be regarded as a
motivic counterpart of the Segal theorem. Also, an alternative
approach to the classical Morel--Voevodsky~\cite{MV} stable homotopy theory $SH(k)$
is suggested in~\cite{GP2a}, which is based on the machinery of framed bispectra.
One of the key
steps in the computations of~\cite{GP1,GP2a} is Theorem~A
proved in this paper. Theorem~A is the main result of the present
paper. In order to state it, we have to recall some definitions and
constructions from~\cite{GP1}.

The framed motive of $X\in Sm/k$ is an
explicitly constructed $S^1$-spectrum $M_{fr}(X)$, which is
connected and an $\Omega$-spectrum in positive degrees
(see~\cite{GP1} for details).
Following the notation of~\cite[Section~8]{GP1}
let $\bb G$ be the cone $(\bb G_m)_+//pt_+$ of the embedding
$pt_+\bl 1\hookrightarrow(\bb G_m)_+$ in the category of pointed
simplicial presheaves $sPre_\bullet(Sm/k)$.
Its sheafification is represented in the category $\Delta^{\op}(Fr_0(k))$ by the object
$\bb G_m^{\wedge 1}$ (see~\cite[Notation 8.1]{GP1}).
For any integer $n\geq 1$ let $\bb G_m^{\wedge n}$ be the $n$th monoidal power of $\bb G_m^{\wedge 1}$
in the symmetric monoidal category $\Delta^{\op}(Fr_0(k))$
(see~\cite[Notation 8.1]{GP1}).
For a variety $X\in Sm/k$ let $M_{fr}(X\times \bb G_m^{\wedge n})$ be the framed motive of the simplicial object
$X\times \bb G_m^{\wedge n}\in \Delta^{\op}(Fr_0(k))$.
It is an explicitly constructed $S^1$-spectrum which is
connected and an $\Omega$-spectrum in positive degrees
(see~\cite[Sections~5 and~6]{GP1} for details).
For brevity we also write $M_{fr}(X)(n)$ to denote $M_{fr}(X\times \bb G_m^{\wedge n})$ and call
$M_{fr}(X)(n)$ the {\it $n$-twisted framed
motive of $X$\/} (see~\cite[Section~11]{GP1}).
The main object of~\cite{GP1} is the bispectrum
   $$M_{fr}^{\bb G}(X)=(M_{fr}(X),M_{fr}(X)(1),M_{fr}(X)(2),\ldots),$$
each term of which is a twisted framed motive
of $X$ and structure maps of the bispectra
   $$M_{fr}(X)(n)\to\underline{\Hom}(\mathbb G,M_{fr}(X)(n+1)),$$
are defined in~\cite[Section~11]{GP1} (they use~\cite[``General
Framework" of Section~5]{GP1}).

The major property of the bispectrum $M_{fr}^{\bb G}(X)$ is
that its levelwise Nisnevich local stable replacement $M^{\bb G}_{fr}(X)_f$
is a fibrant replacement of the suspension bispectrum
$\Sigma_{\bb G}^\infty\Sigma_{S^1}^\infty X_+$. This may also be
viewed as a motivic version of the Barratt, Priddy, and Quillen theorem.
The proof of this major property is given in [GP1] and is heavily based on
the Cancellation Theorem.

The main purpose of the paper is to prove the following (cf.
Voevodsky~\cite{V1})

\begin{theorema}[Cancellation]\label{cancel}
Let $k$ be an infinite perfect field, $X\in Sm/k$ and $n\geq 0$.
Then the following statements are true:

\begin{enumerate}
\item the natural map of $S^1$-spectra
   $$M_{fr}(X)(n)\to\underline{\Hom}(\mathbb G,M_{fr}(X)(n+1))$$
is a schemewise stable equivalence;

\item the induced map of $S^1$-spectra
   $$M_{fr}(X)(n)_f\to\underline{\Hom}(\mathbb G,M_{fr}(X)(n+1)_f)$$
is a schemewise stable equivalence. Here
$M_{fr}(X)(n)_f$ and $M_{fr}(X)(n+1)_f$ are Nisnevich local stable fibrant
replacements of $M_{fr}(X)(n)$ and $M_{fr}(X)(n+1)$ in the injective
local stable model structure of $S^1$-spectra.
\end{enumerate}
\end{theorema}

As an application of Theorem~A we prove the following

\begin{theoremb}\label{bispectrum}
Let $k$ be an infinite perfect field, $X\in Sm/k$ and $n\geq 0$.
Then the bispectrum
   $$M_{fr}^{\bb G}(X)_f=(M_{fr}(X)_f,M_{fr}(X)(1)_f,M_{fr}(X)(2)_f,\ldots)$$
obtained from $M_{fr}^{\bb G}(X)$ by taking levelwise
Nisnevich local stable fibrant replacements with structure maps those of
Theorem~A(2) is a motivically fibrant $(S^1,\bb G)$-bispectrum.
\end{theoremb}

The main strategy of proving Theorem~A is to reduce it to the
``Linear Cancellation Theorem". In order to formulate it, recall
from~\cite[Definition~8.3]{GP1} that the category $\ZF_*(k)$ is the
additive category whose objects are those of $Sm/k$ with Hom-groups
described in Definition~\ref{stab}. Briefly speaking, for every $n\geq 0$ and $X,Y\in Sm/k$
we set
   $$\ZF_n(X,Y):=\bb Z\Fr_n(X,Y)/\langle Z_1\sqcup Z_2-Z_1-Z_2\rangle,$$
where $Z_1,Z_2$ are supports of framed correspondences level $n$
in the sense of Voevodsky~\cite{V2} (see Definition~\ref{stab} as
well). In other words, $\ZF_n(X,Y)$ is the free abelian group
generated by the framed correspondences of level $n$ with connected
supports. We then set
   $$\Hom_{\ZF_*(k)}(X,Y):=\bigoplus_{n\geq 0}\ZF_n(X,Y).$$

Given smooth varieties $X,Y\in Sm/k$ and $n\geq 0$, there is a
canonical suspension morphism $\Sigma:\ZF_n(X,Y)\to\ZF_{n+1}(X,Y)$.
We can stabilise in the $\Sigma$-direction to get an abelian group
(see Definition~\ref{d:stable_linear_fr_corr})
   $$\ZF(X,Y):=\colim(\ZF_0(X,Y)\xrightarrow{\Sigma}\ZF_1(X,Y)\xrightarrow{\Sigma}\cdots).$$
The presheaf $\ZF(Y):=\ZF(-,Y)$ has a canonical structure of a $\ZF_*(k)$-presheaf.
For each scheme $Y\in Sm/k$ and each scheme $S\in Sm/k$ pointed at a $k$-rational point $s\in S$,
the natural functor
   $$\boxtimes: Pre_{Ab}(\ZF_*(k)) \times Pre_{Ab}(\ZF_0(k))\to Pre_{Ab}(\ZF_*(k))$$
defined on p.~\pageref{gfr}
takes the pair $(\ZF(Y),(S,s))$ to the $\ZF_*(k)$-presheaf $\ZF(Y)\boxtimes (S,s)$ which we
also denote by
$\ZF(Y\wedge (S,s))$. By the General Framework of
\cite[Section 5]{GP1} (also see p.~\pageref{gfr})
one has a $\ZF_*(k)$-presheaf
$\underline{\Hom}((S,s),\ZF(Y\wedge (S,s)))$
together with a morphism of $\ZF_*(k)$-presheaves
$$\ZF(Y)\xrightarrow{-\boxtimes (S,s)} \underline{\Hom}((S,s),\ZF(Y\wedge (S,s)).$$

The Linear Cancellation Theorem is formulated as follows
(see Section~\ref{preliminaries} for details).

\begin{theoremc}[Linear Cancellation]
\label{Main}
Let $k$ be an infinite perfect field and let
$Y$ be a $k$-smooth scheme. Then
\[
-\boxtimes (\Gm,1) \colon \ZF(\Delta^\bullet \times -,Y) \to
\underline{\Hom}((\Gm,1),\ZF(\Delta^\bullet \times -,Y\wedge (\Gm,1))
\]
is a quasi-isomorphism of complexes of $\ZF_*(k)$-presheaves of abelian groups.
Here $\Delta^\bullet$ is the standard cosimplicial object in $Sm/k$.
\end{theoremc}




One of the main computational results of~\cite{GNP} says that
schemewise homology of the complex $\ZF(\Delta^\bullet \times-,Y)$
computes homology of the framed motive $M_{fr}(Y)$ of $Y\in Sm/k$.
Moreover, the complex represents the ``linear framed motive" of $Y$
(see~\cite{GNP} for details).


Throughout the paper the base field $k$ is supposed to be infinite. We also employ the following notation:
\begin{itemize}
\item[$\bullet$] all schemes are separated Noetherian $k$-schemes, all morphisms of schemes are $k$-morphisms;
write $\pt$ for the scheme $\textrm{Spec}(k)$.
\item[$\bullet$] $Sm/k$ is the category of smooth $k$-schemes of finite type;
\item[$\bullet$] we refer to the objects of $Sm/k$ as $k$-smooth schemes or smooth $k$-schemes;
\item[$\bullet$] 
Following~\cite{GrD},
by an essentially smooth $k$-scheme we mean a Noetherian $k$-scheme $X$ which is the
inverse limit of a left filtering system $(X_i)_{i\in I}$ with each transition morphism
$X_i \to X_j$ being an \'{e}tale affine morphism between smooth $k$-schemes.
\end{itemize}

\section{Preliminaries}\label{preliminaries}
In this section we collect basic facts for framed correspondences.
We start with preparations.

Let $V$ be a scheme and $Z$ be a closed subscheme. Recall that an
{\it \'{e}tale neighborhood of $Z$ in $V$\/} is a triple
$(W',\pi':W'\to V,s': Z\to W')$ satisfying the following conditions:

\indent (i) $\pi'$ is an \'{e}tale morphism;

\indent (ii) $\pi'\circ s'$ coincides with the inclusion
$Z\hookrightarrow V$ (thus $s'$ is a closed embedding);

\indent (iii) $(\pi')^{-1}(Z)=s'(Z)$.

A morphism between two \'{e}tale neighborhoods
$(W',\pi',s')\to(W'',\pi'',s'')$ of $Z$ in $V$
is a morphism $\rho:W'\to W''$ such that $\pi''\circ\rho=\pi'$ and
$\rho\circ s'=s''$. Note that such $\rho$ is automatically \'etale
by~\cite[VI.4.7]{LNM146}.

\begin{definition}[Voevodsky~\cite{V2}]\label{fr_corr}{\rm


For $k$-smooth schemes $X, Y$ and $n\geq 0$ an {\it explicit framed
correspondence\/}  $\Phi$ of level $n$ consists of the following
data:

\begin{enumerate}
\item a closed subset $Z$ in $\mathbb A^n_X$ which is finite over $X$;
\item an etale neighborhood $p:U\to\mathbb A^n_X$ of $Z$ in $\mathbb A^n_X$;
\item a collection of regular functions $\phi=(\phi_1,\ldots,\phi_n)$ on $U$
such that $\cap_{i=1}^n \{\phi_i=0\}=Z$;
\item a morphism $g:U\to Y$.
\end{enumerate}
The subset $Z$ will be referred to as the {\it support\/} of the
correspondence. We shall also write triples $\Phi=(U,\phi,g)$ or
quadruples $\Phi=(Z,U,\phi,g)$ to denote explicit framed
correspondences.

Two explicit framed correspondences $\Phi$ and $\Phi'$ of level $n$
are said to be {\it equivalent\/} if they have the same support and
there exists an etale neighborhood $V$ of $Z$ in $U \times_{\bb
A^n_X} U'$ such that the morphism $g\circ pr$ agrees with $g' \circ
pr'$ and $\phi \circ pr$ agrees with $\phi' \circ pr'$ on $V$. A
{\it framed correspondence of level $n$} is an equivalence class of
explicit framed correspondences of level $n$.

}\end{definition}

We let $\Fr_n(X,Y)$ denote the set of framed correspondences from
$X$ to $Y$. It as a pointed set with the distinguished point being
the class $0_n$ of the explicit correspondence with $U=\emptyset$.

As an example, the set $\Fr_0(X,Y)$ coincides with the set of
pointed morphisms $X_+\to Y_+$. In particular, for a connected
scheme $X$ one has\label{fr0}
   $$\Fr_0(X,Y)=\text{Hom}_{Sm/k}(X,Y)\sqcup\{0_0\}.$$

If $f:X'\to X$ is a morphism of schemes and $\Phi=(U,\phi,g)$ an
explicit correspondence from $X$ to $Y$ then
   $$f^*(\Phi):=(U'=U\times_X X',\phi\circ pr,g\circ pr)$$
is an explicit correspondence from $X'$ to $Y$.


The following definition is to describe compositions of framed
correspondences.

\begin{definition}\label{Defn_Fr_*}{\rm
Let $X,Y$ and $S$ be $k$-smooth schemes and let
$a=(Z,U,(\varphi_1,\varphi_2,\dots,\varphi_n),g)$
be an explicit correspondence of level $n$ from $X$
to $Y$ and let
$b=(Z',U',(\psi_1,\psi_2,\dots,\psi_m),g')$
be an explicit correspondence of level $m$ from $Y$ to $S$. We
define their composition as an explicit correspondence of level
$n+m$ from $X$ to $S$ by
\[
(Z\times_Y Z',U\times_Y U',(\varphi_1,\varphi_2,\dots,\varphi_n,\psi_1,\psi_2,\dots,\psi_m),g').
\]
Clearly, the composition of explicit correspondences respects the
equivalence relation on them and defines associative pairings
   \begin{equation*}\label{compos}
    \Fr_n(X,Y)\times \Fr_m(Y,S)\to \Fr_{n+m}(X,S).
   \end{equation*}

}\end{definition}

Given $X, Y\in Sm/k$, denote by $\Fr_*(X,Y)$ the set $\bigsqcup_n
\Fr_n(X,Y)$. The composition of framed correspondences defined above
gives a category $\Fr_*(k)$. Its objects are those of $Sm/k$ and the
morphisms are given by the sets $\Fr_*(X,Y)$, $X, Y\in Sm/k$. Since
the naive morphisms of schemes can be identified with certain framed
correspondences of level zero, we get a canonical functor
   $$Sm/k\to \Fr_*(k).$$
One can easily see that for a framed correspondence $\Phi:X\to Y$
and a morphism $f:X'\to X$, one has $f^*(\Phi)=\Phi\circ f$.

\begin{definition}\label{d:boxpairing}{\rm
Let $X,Y,S$ and $T$ be smooth schemes. There is an \textit{external
product}
\[
\Fr_n(X,Y)\times \Fr_m(S,T) \xrightarrow{-\boxtimes -}
\Fr_{n+m}(X\times S, Y\times T)
\]
given by
\begin{gather*}
(Z,U,(\varphi_1,\varphi_2,\dots,\varphi_n),g)\boxtimes
(Z',U',(\psi_1,\psi_2,\dots,\psi_m),g')=
\\
(Z\times Z',U\times U',(\varphi_1,\varphi_2,\dots,\varphi_n,\psi_1,\psi_2,\dots,\psi_m),g\times g').
\end{gather*}

For the constant morphism $c\colon \A^1\to \pt$, we set (following
Voevodsky~\cite{V2})
\[
\Sigma=-\boxtimes (t,c,\{0\},\A^1,t,c)\colon \Fr_n(X,Y)\to \Fr_{n+1}(X,Y)
\]
and refer to it as the \textit{suspension}.
If there is no likelihood of confusion,
we shall also write $\Sigma$ to denote the element $1\cdot(t,c,\{0\},\A^1,t,c)$
in $\ZF_1(\pt,\pt)$ and $\Sigma^n$ for $\Sigma\boxtimes ... \boxtimes \Sigma$
in $\ZF_n(\pt,\pt)$. It will always be clear from the context which of the meanings
for $\Sigma$ is used (either as the suspension or as the element in $\ZF_1(\pt,\pt)$).

Also, following Voevodsky~\cite{V2}, one puts
\[
\Fr(X,Y)=\colim(\Fr_0(X,Y)\xrightarrow{\Sigma}\Fr_1(X,Y)\xrightarrow{\Sigma}\dots
\xrightarrow{\Sigma} \Fr_n(X,Y)\xrightarrow{\Sigma}\dots)
\]
and refer to it as the \textit{set of stable framed correspondences}.
The above external product induces external products
\begin{gather*}
\Fr_n(X,Y)\times \Fr(S,T) \xrightarrow{-\boxtimes -} \Fr(X\times S, Y\times T),\\
\Fr(X,Y)\times \Fr_0(S,T) \xrightarrow{-\boxtimes -} \Fr(X\times S, Y\times T).
\end{gather*}

}\end{definition}

Recall now the definition of the {\it category of linear framed
correspondences $\ZF_*(k)$}.

\begin{definition}\label{stab}{\rm
(see \cite{GP1}) Let $X$ and $Y$ be smooth schemes. Denote by
\begin{itemize}
\item[$\diamond$]
$\ZFr_n(X,Y):=\widetilde{\mathbb{Z}}[\Fr_n(X,Y)]=\mathbb{Z}[\Fr_n(X,Y)]/\mathbb{Z}\cdot
0_n$, i.e the free abelian group generated by the set $\Fr_n(X,Y)$
modulo $\mathbb{Z}\cdot 0_n$;

\item[$\diamond$]
$\ZF_n(X,Y):=\ZFr_n(X,Y)/A$, where $A$ is a subgroup generated by
the elementts
\begin{multline*}
(Z\sqcup Z', U,(\varphi_1,\varphi_2,\dots,\varphi_n),g) - \\
-(Z, U\setminus
Z',(\varphi_1,\varphi_2,\dots,\varphi_n)|_{U\setminus
Z'},g|_{U\setminus Z'}) - (Z',{U\setminus
Z},(\varphi_1,\varphi_2,\dots,\varphi_n)|_{U\setminus
Z},g|_{U\setminus Z}).
\end{multline*}
\end{itemize}
We shall also refer to the latter relation as the {\it additivity
property for supports}. In other words, it says that a framed
correspondence in $\ZF_n(X,Y)$ whose support is a disjoint union
$Z\sqcup Z'$ equals the sum of the framed correspondences with
supports $Z$ and $Z'$ respectively. Note that $\ZF_n(X,Y)$ is
$\mathbb{Z}[\Fr_n(X,Y)]$ modulo the subgroup generated by the
elements as above, because $0_n=0_n+0_n$ in this quotient group,
hence $0_n$ equals zero. Indeed, it is enough to observe that the
support of $0_n$ equals $\emptyset\sqcup\emptyset$ and then apply
the above relation to this support.

The elements of $\ZF_n(X,Y)$ are called {\it linear framed
correspondences of level $n$} or just {\it linear framed
correspondences}.

Denote by $\ZF_*(k)$ the additive category whose objects are those of
$Sm/k$ with $\text{Hom}$-groups defined as
   $$\text{Hom}_{\ZF_*(k)}(X,Y)=\bigoplus_{n\geq 0}\ZF_n(X,Y).$$
The composition is induced by the composition in the category
$\Fr_*(k)$.
Denote by $Pre_{Ab}(\ZF_*(k))$ the Grothendieck category of additive presheaves of abelian groups
on the additive category $\ZF_*(k)$.

Denote by $\ZF_0(k)$ the additive category whose objects are those of
$Sm/k$ with $\text{Hom}$-groups defined as
$\text{Hom}_{\ZF_0(k)}(X,Y)=\ZF_0(X,Y)$.
Clearly, $\ZF_0(k)$ is an additive subcategory of the additive category $\ZF_*(k)$.
Finally, denote by $Pre_{Ab}(\ZF_0(k))$ the category of additive presheaves of abelian groups
on the additive category $\ZF_0(k)$.

There is a natural functor from $Sm/k$ to $\ZF_0(k)$. It is the identity on
objects and takes a regular morphism $f: X\to Y$ to the linear
framed correspondence $1\cdot(X,X\times \A^0,pr_{\A^0},f\circ pr_X)
\in \ZF_0(k)$.

}\end{definition}

\begin{definition}\label{d:boxpairing_linear}{\rm
Let $X,Y,S$ and $T$ be schemes. The external product from
Definition~\ref{d:boxpairing} induces a unique external product
\[
\ZF_n(X,Y)\times \ZF_m(S,T) \xrightarrow{-\boxtimes -}
\ZF_{n+m}(X\times S, Y\times T)
\]
such that for any elements $a \in \Fr_n(X,Y)$ and $b\in \Fr_m(S,T)$
one has $1\cdot a\boxtimes 1\cdot b=1\cdot(a\boxtimes b)\in
\ZF_{n+m}(X\times S, Y\times T)$.

}\end{definition}

\begin{definition}\label{d:stable_linear_fr_corr}{\rm
For any $k$-smooth variety $Y$, the presheaf represented by $Y$
is denoted by $\ZF_*(-,Y)$. One of the main
$\ZF_*(k)$-presheaves of this paper is defined as
   $$\ZF(-,Y)=\colim(\ZF_0(-,Y)\xrightarrow{\Sigma}\ZF_1(-,Y)\xrightarrow{\Sigma}\dots \xrightarrow{\Sigma} \ZF_n(-,Y)\xrightarrow{\Sigma}\dots).$$
For a $k$-smooth variety $X$, the elements of $\ZF(X,Y)$ are also called
{\it stable linear framed correspondences}. Notice that stable linear framed
correspondences do not form morphisms of a category.

}\end{definition}

\noindent{\bf General Framework}.
The pairing $\boxtimes$ of Definition \ref{d:boxpairing_linear} gives rise to a functor
   $$\ZF_*(k) \times \ZF_0(k)\xrightarrow{\boxtimes} \ZF_*(k)$$
taking a pair of schemes $(X,S)$ to $X\times S$ and taking a pair of morphisms
$(a,b)$ to the morphism $a\boxtimes b$. It is naturally extended to a functor\label{gfr}
$$Pre_{Ab}(\ZF_*(k)) \times Pre_{Ab}(\ZF_0(k))\xrightarrow{\boxtimes} Pre_{Ab}(\ZF_*(k)).$$
Given schemes $Y,S\in \ZF_0(k)$, consider the presheaf $\ZF(S)$ in $Pre_{Ab}(\ZF_0(k))$
and the presheaf $\ZF(Y)$ in $Pre_{Ab}(\ZF_*(k))$.
Similarly to~\cite[General Framework, Section 5]{GP1}
there are defined $\ZF_*(k)$-presheaves
$\ZF(Y)\boxtimes S$ and $\underline{\Hom}(S,\ZF(Y)\boxtimes S)$ as well as a natural $\ZF_*(k)$-morphism
$\ZF(Y)\xrightarrow{-\boxtimes S} \underline{\Hom}(S,\ZF(Y)\boxtimes S)$.
By construction, $\ZF(Y)\boxtimes S=\ZF(Y\times S)$. Thus one has the following morphism of $\ZF_*(k)$-presheaves
   $$-\boxtimes id_S: \ZF(Y)\xrightarrow{} \underline{\Hom}(S,\ZF(Y\times S))$$
taking $a\in \ZF(X,Y)$ to $a\boxtimes id_S\in \ZF(X\times S,Y\times S)$.

\begin{definition}\label{d:boxtimes}{\rm
Let $(S,s)$ be a $k$-smooth pointed scheme. Then the morphism $e_s\colon S\xrightarrow{} \pt \xrightarrow{s} S$ defines an
idempotent $\underline{\Hom}(S,e_s): \underline{\Hom}(S,\ZF(Y\times S))\to \underline{\Hom}(S,\ZF(Y\times S))$
in the category of $\ZF_*$-presheaves. Set,
$$\underline{\Hom}(S,\ZF(Y\wedge (S,s))):= \textrm{Ker}[\underline{\Hom}(S,e_s)].$$
Consider the idempotent $\underline{\Hom}(e_s,\underline{\Hom}(S,\ZF(Y\wedge (S,s))))$
of $\underline{\Hom}(S,\ZF(Y\wedge (S,s)))$
in the category of $\ZF_*(k)$-presheaves. Set,
   $$\underline{\Hom}((S,s),\ZF(Y\wedge (S,s))):=\textrm{Ker}[\underline{\Hom}(e_s,\underline{\Hom}(S,\ZF(Y\wedge (S,s))))].$$
For any $X\in Sm/k$ denote by $\ZF(X\wedge (S,s),Y\wedge (S,s))$ the value of $\underline{\Hom}((S,s),\ZF(Y\wedge (S,s)))$ on $X$.
There is a natural morphism of $\ZF_*(k)$-presheaves
   $$-\boxtimes id_{(S,s)}: \ZF(Y)\xrightarrow{} \underline{\Hom}((S,s),\ZF(Y\wedge (S,s))).$$
}\end{definition}


\begin{definition}\label{idemp}{\rm
Let $X$ and $Y$ be $k$-smooth schemes and let $(S,s)$
be a $k$-smooth pointed scheme.
\begin{itemize}
\item[$\diamond$]
Denote by $e_s\colon S\xrightarrow{} \pt \xrightarrow{s} S$ the
idempotent in $\text{End}_{\ZF_0(k)}(S)=\ZF_0(S,S)$
given by the composition of the constant map and the
embedding of $s$ into $S$.

\item[$\diamond$]
For each integer $m\geq 0$ define $\ZF_m(X\wedge (S,s),Y \wedge (S,s))$ as a subgroup of the group
$\ZF_m(X\times S, Y\times S)$ consisting of all $a$ such that $a\circ
(\id_X\boxtimes e_s)=(\id_Y\boxtimes e_s)\circ a=0$.
Note that the suspension map $\Sigma: \ZF_m(X\times S, Y\times S)\to \ZF_{m+1}(X\times S, Y\times S)$
takes the subgroup $\ZF_m(X\wedge (S,s), Y\wedge (S,s)))$ to the subgroup
$\ZF_{m+1}(X\wedge (S,s), Y\wedge (S,s))$. Set,
\begin{multline*}
  \ZF(X\wedge (S,s),Y\wedge (S,s)):=\colim[\ZF_0(X\wedge (S,s),Y\wedge (S,s))\xrightarrow{\Sigma}\ZF_1(X\wedge (S,s),Y\wedge (S,s))\xrightarrow{\Sigma}\dots ],
\end{multline*}
\end{itemize}
It is easy to see that the morphisms
$id_X\boxtimes (\id_\Gm -~ e_s): \ZF_m(X,Y)\to \ZF_m(X\times S,Y\times S)$
take values in $\ZF_m(X\wedge (S,s),Y\wedge (S,s))$. They are compatible with the 
suspension $\Sigma$ and we define a morphism
$$id_X\boxtimes \id_{(S,s)} \colon \ZF(X,Y) \xrightarrow{} \ZF(X\wedge (S,s),Y\wedge (S,s)).$$
}
\end{definition}

\begin{lemma}\label{two_maps_coincides}
Let $Y$ be $k$-smooth scheme and let $(S,s)$
be a $k$-smooth pointed scheme. Then
one has a commutative diagram of $\ZF_*(k)$-presheaves
$$\xymatrix{\ZF(-,Y) \ar[rr]^(.35){-\boxtimes (\id_\Gm -~ e_s)} \ar[d]_{=} && \ZF(-\wedge (S,s),Y\wedge (S,s))  \ar[d]^{can}\\
              \ZF(-,Y) \ar[rr]^(.33){-\boxtimes \id_{(\Gm,1)}}&&  \underline{\Hom}((S,s),\ZF(-,Y\wedge (S,s))),}$$
where $can$ is the canonical isomorphism.
\end{lemma}


\begin{theoremc}\label{Theorem_C2}
Let $X$ and $Y$ be $k$-smooth schemes and let $(\Gm,1)$ be the scheme $\Gm$ pointed by the point $1$.  Then
the morphisms
\begin{equation}\label{first_form}
-\boxtimes (\id_\Gm - e_1) \colon \ZF(\Delta^\bullet \times -,Y) \to
\ZF((\Delta^\bullet \times -)\wedge (\Gm,1),Y\wedge (\Gm,1))
\end{equation}
\begin{equation}\label{second_form}
-\boxtimes \id_{(\Gm,1)} \colon \ZF(\Delta^\bullet \times -,Y) \xrightarrow{}
\underline{\Hom}((\Gm,1),\ZF(\Delta^\bullet \times -,Y\wedge (\Gm,1)))
\end{equation}
are sectionwise quasi-isomorphisms of complexes of $\ZF_*(k)$-presheaves of abelian groups.
\end{theoremc}

\begin{remark}{\rm
By Lemma \ref{two_maps_coincides} the morphism \eqref{first_form} is a quasi-isomorphism
if and only if the morphism \eqref{second_form} is a quasi-isomorphism. Sometimes it is convenient
to work with the morphism \eqref{first_form} and sometimes it is convenient
to work with the morphism \eqref{second_form}.
}
\end{remark}

Theorem~C will be proved at the end of the paper.

\section{Theorem~A and Theorem~B}

Before proving Theorem~A we recall some definitions and
constructions for framed motives for the convenience of the reader. We adhere to~\cite{GP1}.
Let $\Fr_0(k)$ be the category whose objects are those of $Sm/k$ and
whose morphism set between $X$ and $Y$ is given by
the set of framed correspondences of level zero \cite[Example 2.1]{V2}, \cite[Definition 2.1]{GP1}.
As it is shown in~\cite[Section~5]{GP1}, the category of framed correspondences
of level zero $\Fr_0(k)$ has an action by finite pointed sets $Y\otimes
K:=\bigsqcup_{K\setminus *}Y$ with $Y\in Sm/k$ and $K$ a finite pointed
set. Let $U,X\in \Fr_0(k)$.
By the
Additivity Theorem of~\cite{GP1}
the $\Gamma$-space in the sense of
Segal~\cite{Seg}
   $$K\in\Gamma^{\textrm{op}}\mapsto C_*\Fr(U,X\otimes K):=\Fr(U\times\Delta^\bullet,X\otimes K)$$
is special.

\begin{definition}[see~\cite{GP1}]\label{frmotive}{\rm
The {\it framed motive $M_{fr}(X)$ of a smooth algebraic variety\/}
$X\in \Fr_0(k)$ is the Segal $S^1$-spectrum
$(C_*\Fr(-,X),C_*\Fr(-,X\otimes S^1),C_*\Fr(-,X\otimes S^2),\ldots)$
associated with the special $\Gamma$-space
$K\in\Gamma^{\textrm{op}}\mapsto C_*\Fr(-,X\otimes K)$. The framed
motive $M_{fr}(X)\in Sp_{S^1}(k)$ is covariantly functorial in framed
correspondences of level zero.

}\end{definition}


Let $\Fr_0(k)$ be the category whose objects are those of $Sm/k$ and
whose morphism set between $X$ and $Y$ is given by
the set of framed correspondences of level zero \cite[Example 2.1]{V2}, \cite[Definition 2.1]{GP1}.
As it is shown in~\cite[Section~5]{GP1}, the category of framed correspondences
of level zero $\Fr_0(k)$ has an action by finite pointed sets $Y\otimes
K:=\bigsqcup_{K\setminus *}Y$ with $Y\in Sm/k$ and $K$ a finite pointed
set. The cone of $Y$ is the simplicial object $Y\otimes I$ in
$\Fr_0(k)$, where $(I,1)$ is the pointed simplicial set $\Delta[1]$
with basepoint 1. There is a natural morphism $i_0:Y\to Y\otimes I$
in $\Delta^{\textrm{\op}}\Fr_0(k)$. Given a closed inclusion of smooth
schemes $j: Y\hookrightarrow X$, denote by $X//Y$ a simplicial
object in $\Fr_0(k)$ which is obtained by taking the pushout of the diagram
   $X\hookleftarrow Y\bl{i_0}\hookrightarrow Y\otimes I$
in $\Delta^{\textrm{\op}}\Fr_0(k)$.
The simplicial object $X//Y$ termwise equals
   $X,X\sqcup Y,X\sqcup Y\sqcup Y,\ldots$ .
By
$\bb G_m^{\wedge 1}$ we mean the simplicial object $\bb G_m//\{1\}$ in $\Fr_0(k)$.
It looks termwise as
   $$\bb G_m,\bb G_m\sqcup pt,\bb G_m\sqcup pt\sqcup pt,\ldots$$
Applying $M_{fr}(X\times-)$ to $\bb G_m^{\wedge 1}$ and realizing by
taking diagonals, one gets a framed $S^1$-spectrum
$M_{fr}(X\times\bb G_m^{\wedge 1})$. We shall also denote it by
$M_{fr}(X)(1)$. The $n$th iteration gives the spectrum
$M_{fr}(X\times\bb G_m^{\wedge n})$, also denoted by $M_{fr}(X)(n)$.

Similarly to the General Framework on p.~\pageref{gfr} there is a natural pairing
   $$\boxtimes:sPre_\bullet^{fr}(k)\times sPre_\bullet(\Fr_0(k))\to sPre_\bullet^{fr}(k),$$
where $sPre_\bullet^{fr}(k)$ (respectively $sPre_\bullet(\Fr_0(k))$) is the category of pointed simplicial presheaves with framed correspondences
(respectivley the category of pointed simplicial presheaves on $\Fr_0(k)$).
It is extended from the pairing $\Fr_*(k)\times \Fr_0(k) \xrightarrow{\boxtimes} \Fr_*(k)$
that takes $(X,Y)$ to $X\times Y$ and $a\in \Fr_m(X,X')$, $b\in \Fr_0(Y,Y')$ to
$a\boxtimes b\in \Fr_m(X\times X',Y\times Y')$.

We will also write $\wedge$ for the monoidal product in $\Fr_0(k)$ and in $\Delta^{\op}\Fr_0(k)$.
The Yoneda embedding identifies $\Delta^{\op}\Fr_0(k)$ with a full subcategory of $sPre_\bullet(\Fr_0(k))$.
For each integer $n\geq 0$ there is a natural map of spectra
$$a_n: M_{fr}(X\times \bb G_m^{\wedge n})\xrightarrow{- \boxtimes \bb G_m^{\wedge 1}}
\underline{\Hom}(\bb G_m^{\wedge 1},M_{fr}(X\times\bb G_m^{\wedge n+1}))\to
\underline{\Hom}(\bb G,M_{fr}(X\times\bb G_m^{\wedge n+1})),$$
where the right arrow is induced by the adjunction unit
$\textrm{adj}: \bb G\to (\bb G_m^{\wedge 1}|_{Sm/k})$.
Note that $a_n$ respects
framed correspondences of level zero and coincides with the morphism described in~\cite[p.~297]{GP1}).

\begin{definition}\label{d:M_fr_G}
{\rm The $(S^1,\bb G)$-bispectrum $M_{fr}^{\bb G}(X)$ is defined as
$$(M_{fr}(X),M_{fr}(X\times \bb G_m^{\wedge 1} ),M_{fr}(X\times \bb G_m^{\wedge 2}),\ldots)$$
together with the structure morphisms $a_n$-s.
}\end{definition}

We shall prove below (see the proof of Theorem~A) that each $a_n$ is
a schemewise stable equivalence of spectra, but first let us
discuss further useful spectra.
Denote by $\bb ZM_{fr}(X)$, $X\in Sm/k$, the Segal
$S^1$-spectrum
$(\ZFr(\Delta^\bullet\times -,X),\ZFr(\Delta^\bullet\times -,X\otimes S^1),\ldots).$
Denote by $LM_{fr}(X)$ the Segal $S^1$-spectrum
$EM(\ZF(\Delta^\bullet\times-,X))=(\ZF(\Delta^\bullet\times-,X),\ZF(\Delta^\bullet\times-,X\otimes S^1),\ldots).$

The equalities
$\ZF(-,X\sqcup X')=\ZF(-,X)\oplus \ZF(-,X')$ show that the
$\Gamma$-space $(K,*)\mapsto \ZF(\Delta^\bullet\times U,X\otimes K)$ corresponds to the
complex of abelian groups $\ZF(\Delta^\bullet\times U,X)$. Hence $LM_{fr}(X)$ is the
Eilenberg--Mac~Lane spectrum for the complex $\ZF(\Delta^\bullet\times -,X)$.
The $\Gamma$-space
morphism
$$[(K,*)\mapsto \ZFr(\Delta^\bullet\times -,X\otimes K)]\to [(K,*)\mapsto \ZF(\Delta^\bullet\times -,X\otimes K)]$$
induces a morphism of $S^1$-spectra
$l_X:\bb ZM_{fr}(X)\to EM(\ZF(-,X))$.

Note that homotopy groups of
$LM_{fr}(X)=EM(\ZF(\Delta^\bullet\times-,X))$ are equal to homology
groups of the complex $\ZF(\Delta^\bullet\times-,X)$. By~\cite[\S
II.6.2]{Sch} the homotopy groups $\pi_*(\bb ZM_{fr}(X)(U))$ of $\bb
ZM_{fr}(X)$ evaluated at $U\in Sm/k$ are the homology groups
$H_*(M_{fr}(X)(U))$ of $M_{fr}(X)(U)$.

The following result, referred to as the Linearisation Theorem in \cite[Theorem 1.2]{GNP}, is true:

\begin{theorem}[see~\cite{GNP}]\label{cska}
The morphism of $S^1$-spectra
   $$\quad l_X: \bb ZM_{fr}(X)\to LM_{fr}(X)$$
is a schemewise stable equivalence. In particular, if $U$ is
smooth, then
$$H_*(M_{fr}(X)(U))=\pi_*(\bb ZM_{fr}(X)(U))=\pi_*(LM_{fr}(X)(U))=H_*(\ZF(\Delta^\bullet\times U,X)).$$
\end{theorem}

Replacing simplicial framed sheaves $C_*Fr$ in Definition \ref{d:M_fr_G} by
simplicial abelian framed presheaves $C_*\ZF$, we
define Segal $S^1$-spectra $LM_{fr}(X\times\bb G_m^{\wedge n})$-s.
Following the General Framework on p.~\pageref{gfr}, there is a natural 
morphism of $S^1$-spectra  for each integer $n\geq 0$
   \begin{equation}\label{strelki_c}
   c_n: LM_{fr}(X\times\bb G_m^{\wedge n})\xrightarrow{-\boxtimes \bb G_m^{\wedge 1}} \underline{\Hom}(\bb G_m^{\wedge 1},M_{fr}(X\times\bb G_m^{\wedge n+1}))
   \to \underline{\Hom}(\bb G,LM_{fr}(X\times\bb G_m^{\wedge n+1})),
   \end{equation}
where the right arrow is induced by the adjunction unit
$\textrm{adj}: \bb G\to (\bb G_m^{\wedge 1}|_{Sm/k})$.

\begin{definition}\label{d:LM_fr_G}
{\rm The $(S^1,\bb G)$-bispectrum $LM_{fr}^{\bb G}(X)$ is defined as
$$(LM_{fr}(X),LM_{fr}(X\times \bb G_m^{\wedge 1} ),LM_{fr}(X\times \bb G_m^{\wedge 2}),\ldots)$$
together with the structure morphisms $c_n$-s.
}\end{definition}

We are now in a position to prove Theorem~A.

\begin{proof}[Proof of Theorem~A]
(1). We claim that for every $n>0$ the sequence
   $$M_{fr}(X)(n-1)\to M_{fr}(X\times\bb G_m)(n-1)\to M_{fr}(X)(n)$$
is a homotopy cofiber sequence of $S^1$-spectra. Since all spectra are connected,
it is enough to show that
   $$\bb ZM_{fr}(X)(n-1)\to \bb ZM_{fr}(X\times\bb G_m)(n-1)\to \bb ZM_{fr}(X)(n)$$
is a homotopy cofiber sequence of $S^1$-spectra. By Theorem~\ref{cska} the latter
is equivalent to showing that
   $$LM_{fr}(X)(n-1)\to LM_{fr}(X\times\bb G_m)(n-1)\to LM_{fr}(X)(n)$$
is a homotopy cofiber sequence of $S^1$-spectra. This sequence is a homotopy cofiber sequence if and only if
   $$\bb ZF(\Delta^\bullet\times-,X\times\bb G_m^{\wedge(n-1)})\to \bb ZF(\Delta^\bullet\times-,X\times\bb G_m^{\wedge(n-1)}\times\bb G_m)\to
      \bb ZF(\Delta^\bullet\times-,X\times\bb G_m^{\wedge n})$$
is a homotopy cofiber sequence of complexes of abelian presheaves. But this is
obvious because $\bb ZF(\Delta^\bullet\times-,X\times\bb G_m^{\wedge n})$ is the mapping cone
of the left arrow, and hence the desired claim follows. We have used here the fact that
$\bb ZF(-,X\sqcup Y)=\bb ZF(-,X)\oplus\bb ZF(-,Y)$.

Next, it is enough to prove that
   $$a_0:M_{fr}(X)\to\underline{\Hom}(\bb G,M_{fr}(X)(1))$$
is a schemewise equivalence of spectra. Indeed, consider a
commutative diagram of homotopy cofiber sequences in $Sp_{S^1}(k)$
   $$\xymatrix{M_{fr}(X)(n-1)\ar[d]_{a_{n-1}}\ar[r]&M_{fr}(X\times\bb G_m)(n-1)\ar[d]_{a_{n-1}}\ar[r]&M_{fr}(X)(n)\ar[d]^{a_n}\\
               \underline{\Hom}(\bb G,M_{fr}(X)(n))\ar[r]&\underline{\Hom}(\bb G,M_{fr}(X\times\bb G_m)(n))\ar[r]&\underline{\Hom}(\bb
               G,M_{fr}(X)(n+1))}$$
with $n\geq 1$. If $a_{n-1}$ is a schemewise equivalence of spectra,
then so is $a_n$ by~\cite[13.5.10]{Hir}.
Thus using induction in $n$, it suffices to verify that $a_{0}$ is a
schemewise equivalence of spectra.

By the stable Whitehead theorem~\cite[II.6.30]{Sch} $a_0$ is a
stable equivalence if and only if so is
   $$a_0:\bb ZM_{fr}(X)\to\bb Z[\underline{\Hom}(\bb G,M_{fr}(X\times\bb G_m^{\wedge 1}))].$$

Consider a commutative diagram of homotopy fiber sequences in $Sp_{S^1}(k)$
   $$\xymatrix{\underline{\Hom}(\bb G,M_{fr}(X\times\bb G_m^{\wedge 1}))\ar[r]\ar[d]
               &\underline{\Hom}(\bb G_m,M_{fr}(X\times\bb G_m^{\wedge 1}))\ar[d]\ar[r]
               &M_{fr}(X\times\bb G_m^{\wedge 1})\ar[d]\\
               \bb Z[\underline{\Hom}(\bb G,M_{fr}(X\times\bb G_m^{\wedge 1}))]\ar[r]\ar[d]_{\ell_X}
               &\underline{\Hom}(\bb G_m,\bb ZM_{fr}(X\times\bb G_m^{\wedge 1}))\ar[d]\ar[r]
               &\bb ZM_{fr}(X\times\bb G_m^{\wedge 1})\ar[d]^{l_X}\\
               \underline{\Hom}(\bb G,LM_{fr}(X\times\bb G_m^{\wedge 1}))\ar[r]
               &\underline{\Hom}(\bb G_m,LM_{fr}(X\times\bb G_m^{\wedge 1}))\ar[r]
               &LM_{fr}(X\times\bb G_m^{\wedge 1})}$$
The arrow $l_X$ and the middle lower arrow are a stable weak equivalences of
spectra by Theorem~\ref{cska}. It follows that $\ell_X$ is a stable weak equivalence.
Consider a commutative diagram
   $$\xymatrix{\bb ZM_{fr}(X)\ar[d]_{l_X}\ar[r]^(.35){a_0} &\bb Z[\underline{\Hom}(\bb G,M_{fr}(X\times\bb G_m^{\wedge 1}))]\ar[d]^{\ell_X}\\
               LM_{fr}(X)\ar[r]^(.33){c_0}&\underline{\Hom}(\bb G,LM_{fr}(X\times\bb G_m^{\wedge 1})).}$$
Since $l_X,\ell_X$ are stable weak
equivalences, it follows that $a_0$ is a stable local equivalence
if and only if so is $c_0$.
By Theorem D from Appendix A the morphism $c_0$ is a sectionwise stable weak equivalence.
The proof of the first part of the theorem is completed.

(2). Since each spectrum $M_{fr}(X)(n)_f$ is fibrant in the injective local stable model structure of $S^1$-spectra,
it is enough to show that each map
   $$b_n:M_{fr}(X)(n)_f\to\underline{\Hom}(\mathbb G,M_{fr}(X)(n+1)_f),\quad n\geq 0,$$
is a Nisnevich local stable equivalence of spectra. Using the same argument as in the proof of the first statement,
it suffices to verify that $b_0$ is a local stable equivalence.

There is a commutative diagram
   $$\xymatrix{M_{fr}(X)\ar[r]^(.35){a_0}\ar[d]_\alpha&\underline{\Hom}(\mathbb G,M_{fr}(X)(1))\ar[d]^{d_1}\\
                       M_{fr}(X)_f\ar[r]^(.35){b_0}&\underline{\Hom}(\mathbb G,M_{fr}(X)(1)_f)}$$
in which the left vertical arrow is a local stable equivalence and $a_0$ is a schemewise stable equivalence
by the first statement. It follows that $b_0$ is a local stable equivalence if and only if so is $d_1=\underline{\Hom}(\bb G,\alpha)$.

The presheaves of stable homotopy groups of $\underline{\Hom}(\mathbb G,M_{fr}(X)(1))$ equal
$(\pi_n(M_{fr}(X)(1)))_{-1}$. These presheaves are $\bb A^1$-invariant quasi-stable radditive with framed
correspondences (see~\cite[Introduction]{GP2} for the definition of such presheaves).
It follows
from~\cite[Theorem~1.1]{GP2} (complemented by~\cite{DP} in characteristic 2) that
each Nisnevich sheaf $((\pi_n(M_{fr}(X)(1)))_{-1})^{\nis}$ is strictly $\bb A^1$-invariant
quasi-stable radditive with framed correspondences.

Each spectrum $M_{fr}(X)(n)$ has homotopy invariant, quasi-stable radditive
presheaves with framed correspondences of stable homotopy groups $\pi_*(M_{fr}(X)(n))$.
By~\cite[Theorem~1.1]{GP2} (complemented by~\cite{DP} in characteristic 2) the Nisnevich sheaves
$\pi_*^{\nis}(M_{fr}(X)(n))$ are strictly homotopy invariant. It follows from~\cite[Proposition~7.1]{GP1}
that $M_{fr}(X)(n)_f$ is motivically fibrant in the injective stable motivic model structure of $S^1$-spectra.

In order to compute the Nisnevich sheaf $\pi_n^{\nis}(\underline{\Hom}(\mathbb G,M_{fr}(X)(1)_f))$,
consider the Brown--Gersten convergent spectral sequence
   $$H^p_{\nis}(V\times\bb G_m,\pi_q^{\nis}(M_{fr}(X)(1)))\Rightarrow\pi_{q-p}(M_{fr}(X)(1)_f(V\times\bb G_m)),\quad V\in Sm/k.$$
It follows from~\cite[Corollary~16.8, Theorems~17.15-16]{GP2} that each presheaf
$$V\mapsto H^p_{\nis}(V\times\bb G_m,\pi_q^{\nis}(M_{fr}(X)(1)))$$
is $\bb A^1$-invariant quasi-stable radditive with framed correspondences.

Let $V \in Sm/k$ be irreducible, $u\in V$ be a point, $U=\Spec(\cc O_{V,u})$.
Let $U^h_u$ be the henselization of $U$ at $u$ and let
$k(U^h_u)$ be the function field on $U^h_u$.
Consider the above spectral sequence and replace $V$ by $U^h_u$ in it.
We claim that in this case the spectral sequence degenerates and
$H^0_{\nis}(U^h_u\times\bb G_m,\pi_n^{\nis}(M_{fr}(X)(1)))=\pi_{n}(M_{fr}(X)(1)_f(U^h_u\times\bb G_m))$.
For this notice that
by~\cite[3.15(3')]{GP2}
the map
$H^p_{\nis}(\bb G_m\times U,\pi_q^{\nis}(M_{fr}(X)(1)))\hookrightarrow H^p_{\nis}(\bb G_{m,k(U^h_u)},\pi_q^{\nis}(M_{fr}(X)(1)))$
is injective, where $\eta_h:\Spec(k(U^h_u)) \to U^h_u$
is the canonical morphism. In turn,
by \cite[3.15(1)]{GP2} the canonical homomorphism
$$H^p_{\nis}(\bb G_{m,k(U^h_u)},\pi_q^{\nis}(M_{fr}(X)(1)))\hookrightarrow
H^p_{\nis}(\Spec (k(U^h_u)(t)),\pi_q^{\nis}(M_{fr}(X)(1)))$$
is injective. Since
$0=H^p_{\nis}(\Spec (k(U^h_u)(t)),\pi_q^{\nis}(M_{fr}(X)(1)))$
for $p>0$, the group
$H^p_{\nis}(U^h_u\times\bb G_{m},\pi_q^{\nis}(M_{fr}(X)(1)))$ vanishes for $p>0$.
Thus we have checked the equality
$$H^0_{\nis}(U^h_u\times\bb G_m,\pi_n^{\nis}(M_{fr}(X)(1)))=\pi_{n}(M_{fr}(X)(1)_f(U^h_u\times\bb G_m)).$$
We can conclude that
   $\pi_n^{\nis}(\underline{\Hom}(\mathbb G_m,M_{fr}(X)(1)_f))=\pi_n^{\nis}(M_{fr}(X)(1)_f)(\mathbb G_m\times-).$
It follows that
   $$\pi_n^{\nis}(\underline{\Hom}(\mathbb G,M_{fr}(X)(1)_f))=(\pi_n^{\nis}(M_{fr}(X)(1)_f))_{-1}=(\pi_n^{\nis}(M_{fr}(X)(1)))_{-1}.$$

It remains to show that the morphism of $\bb A^1$-invariant radditive quasi-stable framed sheaves
   $$((\pi_n(M_{fr}(X)(1)))_{-1})^{\nis}\to(\pi_n^{\nis}(M_{fr}(X)(1)))_{-1}$$
is an isomorphism. Using~\cite[3.15(3')]{GP2} it suffices to check that it is an isomorphism for every
field extension $K/k$. The homomorphism of abelian groups
   $$((\pi_n(M_{fr}(X)(1)))_{-1})^{\nis}(K)=(\pi_n(M_{fr}(X)(1)))_{-1}(K)\to(\pi_n^{\nis}(M_{fr}(X)(1)))_{-1}(K)$$
is an isomorphism, because for every homotopy invariant radditive quasi-stable framed presheaf of abelian groups
$\cc F$ and every open $V\subset\bb A^1_K$, one has $\cc F(V)=\cc F^{\nis}(V)$ (see the proof of~\cite[3.1]{GP2}).
This completes the proof of Theorem~A.
\end{proof}

\begin{proof}[Proof of Theorem~B]
The proof of Theorem~A(2) shows that $M_{fr}(X)(n)_f$ is motivically fibrant
in the injective stable motivic model structure of $S^1$-spectra.
By Theorem~A each structure
map $b_n$ is a schemewise equivalence. We conclude that the
bispectrum  $M_{fr}^{\bb G}(X)_f$ is a motivically fibrant $(S^1,\bb
G)$-bispectrum in the sense of Jardine~\cite{Jar}.
\end{proof}

\section{Useful lemmas}

In this section we discuss several useful $\mathbb{A}^1$-homotopies
and collect a number of facts used in the following sections. We
start with some definitions and notation.

\begin{definition}{\rm
Let $\mathcal F: Sm/k \to Sets$ be a presheaf of sets. Let $X\in
Sm/k$ be a smooth variety and $a,b \in \mathcal F(X)$ be two
sections. We write $a\sim b$ if $a$ and $b$ are in the same
connected component of the simplicial set $\mathcal
F(\Delta^{\bullet}\times X)$. If $h\in \mathcal F(\Delta^{1}\times
X)$ is such that $\partial_0(h)=a$ and $\partial_1(h)=b$, then we
will write $a\frac{h}{}b$. In this case $a\sim b$.

Let $\mathcal A: Sm/k \to Ab$ be a presheaf of abelian groups. Let
$X\in Sm/k$ be a smooth variety and $a,b \in \mathcal A(X)$ be two
sections. We will write $a\sim b$ if the classes of $a$ and $b$ in
$H_0(\mathcal A(\Delta^{\bullet}\times X))$ coincide. This is
equivalent to saying that there is $h\in \mathcal A(\Delta^{1}\times
X)$ such that $\partial_0(h)=a$ and $\partial_1(h)=b$. For such an
$h$ we will write $a\frac{h}{}b$.

}\end{definition}

\begin{definition}{\rm
Let $\mathcal{F}$ and $\mathcal{G}$ be two presheaves of sets on the
category of $k$-smooth schemes and let $\phi_0,\phi_1: \mathcal{F}
\rightrightarrows \mathcal{G}$ be two morphisms. An {\it
$\A^1$-homotopy\/} between $\phi_0$ and $\phi_1$ is a morphism $H:
\mathcal{F}\to \text{\underline {Hom}}(\A^1,\mathcal{G})$ such that
$H_0=\phi_0$ and $H_1=\phi_1$. We will write $\phi_0\sim \phi_1$ if
there is an $\A^1$-homotopy between $\phi_0$ and $\phi_1$.

Let $\mathcal{A}$ and $\mathcal{B}$ be two presheaves of abelian
groups on the category of $k$-smooth schemes and let $\phi_0,\phi_1:
\mathcal{A} \rightrightarrows \mathcal{B}$ be two morphisms. An {\it
$\A^1$-homotopy\/} between $\phi_0$ and $\phi_1$ is a morphism $H:
\mathcal{A}\to \text{\underline {Hom}}(\A^1,\mathcal{B})$ of
presheaves of abelian groups such that $H_0=\phi_0$ and
$H_1=\phi_1$. If $H$ is an $\A^1$-homotopy between $\phi_0$ and
$\phi_1$, then we will write $\phi_0\frac{H}{} \phi_1$. If we do not
specify an $\A^1$-homotopy between $\phi_0$ and $\phi_1$, then we
will write $\phi_0\sim \phi_1$.

If $\phi: \mathcal{A} \to \mathcal{B}$ is a morphism of presheaves
of abelian groups, then there is a constant $\A^1$-homotopy
$H_{\phi}$ between $\phi$ and $\phi$ defined as follows. Given $a\in
\mathcal A(X)$ set $H_{\phi}(a)=pr^*_X(\phi(a)) \in \mathcal
B(X\times \A^1)$.

}\end{definition}

\begin{lemma} \label{l:homotopyhomology}
Let $\mathcal{A}$ and $\mathcal{B}$ be two presheaves of abelian
groups on the category of $k$-smooth schemes and let $\phi_0,\phi_1:
\mathcal{A} \rightrightarrows \mathcal{B}$ be two morphisms such
that $\phi_0\sim \phi_1$. Then the induced morphisms
$$\phi_0, \phi_1: \mathcal{A}(\Delta^{\bullet})\rightrightarrows \mathcal{B}(\Delta^{\bullet})$$
between two simplicial abelian groups give the same morphisms on the
homology of the associated Moore complexes.
\end{lemma}

\begin{lemma}\label{l:H1_and_H_2}
Let $\phi_0,\phi_1,\phi_2: \mathcal A \to \mathcal B$ be morphisms
of presheaves of abelian groups and let $\phi_0\frac{H'}{}\phi_1$
and $\phi_1\frac{H''}{}\phi_2$. Then
$$\phi_0\frac{H'+H''-H_{\phi_1}}{}\phi_2$$
\end{lemma}

\begin{lemma}
Let $\mathcal{A}$ and $\mathcal{B}$ be two presheaves of abelian
groups on the category of $k$-smooth schemes and let $\phi_0
\frac{H}{}\phi_1$. Let $\rho: \mathcal{A'} \to \mathcal{A}$ be a
morphism. Then $\phi_0\circ \rho \frac{H\circ \rho}{} \phi_1\circ
\rho$. Moreover, let $\eta: \mathcal{B} \to \mathcal{B'}$ be a
morphism, then
$\psi\circ \phi_0 \frac{\psi \circ H}{} \psi \circ
\phi_1$ with $\psi={\underline{\Hom}}(\A^1,\eta):
{\underline{\Hom}}(\A^1,\mathcal B) \to
{\underline{\Hom}}(\A^1,\mathcal B')$.
\end{lemma}

We now want to discuss actions of matrices on framed correspondences
and associated homotopies. Let $X$ and $Y$ be $k$-smooth schemes and
$A\in GL_n(k)$ be a matrix. Then $A$ defines an automorphism
$$\phi_A: \Fr_n(-\times X,Y) \to \Fr_n(-\times X,Y)$$
of the presheaf $\Fr_n(-\times X,Y)$ in the following way. Given
$W\in Sm/k$ and
$a=(Z,U,(\varphi_1,\varphi_2,\dots,\varphi_n),g)
\in \Fr_n(W\times X,Y)$, set
$$
\phi_A (Z,U,(\varphi_1,\varphi_2,\dots,\varphi_n),g))
:= (Z,U,A\circ(\varphi_1,\varphi_2,\dots,\varphi_n),g),
$$
where $A$ is regarded as a linear automorphism of $\A^n_k$.

The automorphism $\phi_A$ of the presheaf $\Fr_n(-\times X,Y)$
induces an automorphism of the free abelian presheaf $\mathbb
Z[\Fr_n(-\times X,Y)]$ and an automorphism $\varphi_A$ of the
presheaf of abelian groups $\ZF_n(-\times X,Y)$.

\begin{definition} \label{d:homotopy_1}{\rm
Let $A\in SL_n(k)$. Choose a matrix $A_s\in SL_n(k[s])$ such that
$A_0=id$ and $A_1=A$. The matrix $A_s$, regarded as a morphism $\A^n
\times \A^1 \to \A^n$, gives rise to an $\A^1$-homotopy $h$ between
$id$ and $\phi_A$ as follows. Given
$a=(\alpha,f,Z,U,\varphi,g))=((\alpha_1,\alpha_2,\dots,\alpha_n),f,Z,U,(\varphi_1,\varphi_2,...,\varphi_n),g)
\in \Fr_n(W\times X,Y)$, one sets
   $$h(a)=(\alpha, f\times id_{\A^1}, Z\times \A^1,U\times \A^1, A_s\circ (\varphi\times id_{\A^1}),g\circ pr_U) \in \Fr_n(W\times X\times \A^1,Y).$$
Clearly, $h_0(a)=a$ and $h_1(a)=\phi_A(a)$. By linearity the
homotopy $h$ induces an $\A^1$-homotopy $H_{A_s}$
$$id\frac{H_{A_s}}{} \varphi_A: \ZF_n(- \times X,Y)\rightrightarrows \ZF_n(- \times X,Y)$$
between the identity $id$ and the morphism $\varphi_A$.

}\end{definition}

\begin{lemma}
Let $\rho: \ZF_m(- \times X,Y) \to \ZF_n(- \times X,Y)$ be a
presheaf morphism. Let $A\in SL_n(k)$, $A_s\in SL_n(k[s])$ and
$H_{A_s}$ be as in Definition~\ref{d:homotopy_1}. Then one has
$$\rho \frac{H_{A_s}\circ \rho}{} \varphi_A \circ \rho: \ZF_m(- \times X,Y) \rightrightarrows \ZF_n(- \times X,Y).$$
\end{lemma}

For $b\in\ZF_m(Y,S)$ define a presheaf morphism
$$\varphi_b: \ZF_n(-\times X,Y) \to \ZF_{n+m}(-\times X,S)$$
sending $a\in \ZF_n(W\times X,Y)$ to $b\circ a\in \ZF_{n+m}(W\times
X,S)$. Also, any $b\in\ZF_m(pt,pt)$ defines a morphism of presheaves
$$-\boxtimes b: \ZF_n(-\times X,Y) \to \ZF_{n+m}(-\times X,Y)$$
sending $a\in \ZF_n(W\times X,Y)$ to $a\boxtimes b\in
\ZF_{n+m}(W\times X,Y)$.

The next three lemmas are straightforward.

\begin{lemma} \label{l:b1_and_b2}
Let $b_1,b_2\in\ZF_m(Y,S)$ be such that $b_1\sim b_2$, then
$$\varphi_{b_1}\sim \varphi_{b_2}: \ZF_n(-\times X,Y) \rightrightarrows \ZF_{n+m}(-\times X,S).$$
\end{lemma}

\begin{lemma} \label{l:box_b1_and_b2}
Let $b_1,b_2\in \ZF_m(pt,pt)$ and $h\in \ZF_m(\A^1,pt)$ be such that
$b_1\frac{h}{} b_2$, then
$$(- \boxtimes {b_1})\frac{-\boxtimes h}{} (- \boxtimes {b_2}): \ZF_n(-\times X,Y) \rightrightarrows \ZF_{n+m}(-\times X,Y).$$
\end{lemma}

The following lemma is proved in Appendix~\ref{henzelization}.

\begin{lemma}\label{Retraction}
Let $z\in \mathbb A^m$ be a $k$-rational point. Set $U'=(\mathbb
A^m)^h_{z}$ to be the henzelization of $\mathbb A^m$ at the point
$z$. Let $i_z: \pt \hookrightarrow U'$ be the closed point of $U'$.
Let $U'_{s}:= (\mathbb A^1\times \mathbb A^m)^h_{\mathbb A^1\times
z}$ be the henzelization of $\mathbb A^1\times \mathbb A^m$ at
$\mathbb A^1\times z$. Then
the morphism $f_s: \mathbb A^1\times \mathbb A^m\to \mathbb A^m$ mapping $(s,y)$
to $s\cdot (y-x)+x$ induces a morphism $H_{s}:=f^h_s: U'_{s}\to U'$ such that:

\begin{itemize}
\item[(a)] $H_1:=(f^h_s)|_{(1\times X)^h_{(1,x)}}: U' \to U'$ is the identity morphism;
\item[(b)] $H_0:=(f^h_s)|_{(0\times X)^h_{(0,x)}}: U' \to U'$ coincides with the composite morphism $U'
\xrightarrow{p^h} \pt \xrightarrow{s_z} U'$, where $p^h: U' \to
\pt=\spec(k)$ is the structure morphism and
$s_z: \pt \hookrightarrow U'$ is the closed point of $U'$.
\end{itemize}
\end{lemma}

Let $z\in \mathbb A^m$ be a $k$-rational point. The projection
$pr: \mathbb A^1\times \mathbb A^m\to \mathbb A^m$
induces a morphism $can_s:=pr^h: U'_{s}\to U'$ such that $can_0=can_1=\id_{U'}$
(see Appendix~\ref{henzelization}). The preceding lemma gives the following

\begin{corollary}\label{Replacing_g}
Let $z\in \mathbb A^m$ be a $k$-rational point and let
$(z,U',\psi;\id_{U'}) \in \Fr_m(\pt,U')$ with $U'$ as in Lemma~\ref{Retraction}.
Suppose $U'_s$ is as in Lemma \ref{Retraction} and let
$h_{s}=(\mathbb A^1\times z,U'_s,can^*_s(\psi);H_{s}) \in \Fr_m(\mathbb
A^1,U')$. Then one has:

\begin{itemize}
\item[(a)] $h_1=(z,U',\psi;\id_{U'}) \in \Fr_m(\pt,U')$;
\item[(b)] $h_0=(z,U',\psi;s_z\circ p^h)=s_z\circ (\{z\},U',\psi;p^h) \in
\Fr_m(\pt,U')$, where $p^h: U' \to \pt=\spec(k)$ is the structure
morphism and $s_z:\pt\hookrightarrow U'$ is the closed point of $U'$.
\end{itemize}
\end{corollary}

\begin{lemma} \label{l:homotopyconstant}
Let $z\in \mathbb A^m$ be a $k$-rational point. Let $Y$ be a
$k$-smooth scheme and let
$(z,U,(\varphi_1,\varphi_2,\hdots,\varphi_m),g)\in \Fr_m(\pt,Y)$ be a
framed correspondence. Then
   $$(z,U,(\varphi_1,\varphi_2,\hdots,\varphi_m),g) \sim
       (z,U,(\varphi_1,\varphi_2,\dots,\varphi_m),c_{g(z)}),$$
where $c_{g(z)}=g(z)\circ p: U\xrightarrow{p} \pt \xrightarrow{g(z)} Y$.
\end{lemma}

\begin{proof}
Let $U'$, $U'_s$, $i_z$ and $h_s$ be as in Corollary~\ref{Replacing_g}.
Let $\pi: U' \to U$ be the canonical morphism.
Set $h'_s=g\circ \pi \circ h_s\in \Fr_m(\pt,Y)$. We want to check
that $h'_1=(z,U,\varphi,g)$ and $h'_0=(z,U,\varphi,c_{g(z)})$. This
will prove our statement. One has,
$$h'_1=(g\circ \pi)\circ h_1=(g\circ \pi)\circ (z,U',\varphi \circ \pi;id_{U'})=(z,U',\varphi\circ \pi;g\circ \pi)=(z,U,\varphi;g),$$
$$h'_0=(g\circ \pi)\circ h_0=(g\circ \pi)\circ (z,U',\varphi \circ \pi;s_z\circ p^h)=(z,U',\varphi \circ \pi;g\circ \pi\circ s_z\circ p^h)=$$
$$=(z,U',\varphi \circ \pi;c_{g(z)}\circ \pi)=(z,U,\varphi;c_{g(z)})$$
as required.
\end{proof}

\begin{lemma} \label{l:canonicalpolynomial}
Let $Y$ be a $k$-smooth scheme and let $(Z,U,\varphi,g)\in
\Fr_1(\pt,Y)$ be a framed correspondence. Suppose that $U\subset
\A^1$ and $\varphi=p(t)\in k[t]$ is a polynomial, where $t$ is the
coordinate function on $\A^1$. Let $g: U\to Y$ be a morphism.
\begin{enumerate}
\item
Then for every $a\in k$ we have
\[
(Z,U,p(t),g(t)) \sim (m_a^{-1}(Z),m_a^{-1}(U),p(t-a),g(t-a)) \in
\Fr_1(\pt,Y),
\]
where $m_a\colon \A^1\to \A^1$ is given by $m_a(t)=t-a$.
\item
If $Z=\{x_0\}$ for some $x_0\in k$ and $p(t)=(t-x_0)^nr(t)$ with $r(t)$ invertible on $U$, then
\[
(Z,U,p(t),g) \sim (\{0\},\A^1,r(x_0)t^{n},c_{g(x_0)}) \in
\Fr_1(\pt,Y),
\]
where $c_{g(x_0)}\colon \A^1\to \pt \xrightarrow{g(x_0)} Y$ is the
constant map taking $\A^1$ to the point $g(x_0)\in Y$.
\end{enumerate}
\end{lemma}

\begin{proof}
(1) The homotopy is given by
\[
(m_{sa}^{-1}(Z),m_{sa}^{-1}(U),p(t-sa),g(t-sa))\in \Fr_1(\A^1,Y),
\]
where $s$ is the homotopy parameter and $m_{sa}\colon
\A^1\times\A^1\to \A^1$ is the morphism $m_{sa}(t)=t-sa$.

(2) Using the preceding statement, we may assume that $x_0=0$.
Consider a polynomial
   $$h(s,t)=sr(t)t^{n}+(1-s)r(0)t^{n}\in k[s,t].$$
If $r_1(t)$ is such that $r(t)=r(0)+ t\cdot r_1(t)$, then one has $h(s,t)=t^n\cdot( r(0)+t\cdot r_1(t)\cdot s)$.
If $S$ is the vanishing locus of
$ r(0)+t\cdot r_1(t)\cdot s$, then $S\cap \A^1\times 0=\emptyset$.
Hence for the zero locus $Z(h)$ of $h$ one has $Z(h)=(\A^1\times 0)\sqcup S$.
The framed correspondence
\[
(\A^1 \times \{0\},(\A^1\times U)\setminus S,
sr(t)t^{n}+(1-s)r(0)t^{n},g\circ pr_U)\in \Fr_1(\A^1,Y)
\]
yields the relation $(\{0\},U,r(t)t^{n},g) \sim
(\{0\},U,r(0)t^{n},g)$ in $\Fr_1(\pt,Y)$.
Lemma~\ref{l:homotopyconstant} shows that
\[
(\{0\},U,r(0)t^{n},g) \sim
(\{0\},U,r(0)t^{n},g(0))=(\{0\},\A^1,r(0)t^{n},g(0)) \in
\Fr_1(\pt,Y)
\]
and our lemma follows.
\end{proof}

\begin{lemma} \label{l:polynomialcycle}
Let $a\in k^{\times}$. Let $p(t),q(t) \in k[t]$ be two polynomials of degree $n$ with the
leading coefficient $a$. Let $(Z(p),\A^1,p(t),c)\in\Fr_1(\pt,\pt)$, $(Z(q),\A^1,q(t),c)\in\Fr_1(\pt,\pt)$
be two framed correspondences. Here $c\colon \A^1\to \pt$ is the structure morphism. Then
   $$(Z(p),\A^1,p(t),c) \sim (Z(q),\A^1,q(t),c) \in \Fr_1(\pt,\pt).$$
\end{lemma}

\begin{proof}
As a polinomial in $t$ the leading coefficient of the polinomial $p(t)+s(q(t)-p(t))$ is $a\in k^{\times}$.
Hence the $k[s]$-module $k[s,t]/(p(t)+s(q(t)-p(t)))$ is a free module rank $n$.
Let $Z_s\subset \A^1\times\A^1$ be the vanishing locus of $p(t)+s(q(t)-p(t))$.
The desired homotopy is given by the framed correspondence
   $$(Z_s,\A^1\times\A^1,p(t)+s(q(t)-p(t)),c'),$$
where $s$ is the homotopy parameter and $c'\colon \A^1\times\A^1\to
\pt$ is the canonical projection.
\end{proof}

\section{Homotopies for swapping coordinates of $\mathbb{G}_m\times\mathbb{G}_m$}

In this section we follow notation of Section~\ref{preliminaries}.
Denote by $\varepsilon=(\{0\},\A^1, -t,c)\in \Fr_1(\pt,\pt)$, where
$c\colon \A^1\to \pt$ is the canonical projection.
We work in this Section with the elements $\Sigma^n \in \ZF_n(\pt,\pt)$
as in Definitition \ref{d:boxpairing}.

\begin{proposition} \label{p:injectivity}
Let $Y$ be a $k$-smooth scheme. Then the canonical homomorphism
\[
H_0(\ZF(\Delta^{\bullet}\times \Gm\times\Gm,Y)) \to
H_0(\ZF(\Delta^{\bullet}_{\Spec k(t,u)},Y))
\]
is injective.
\end{proposition}
\begin{proof}
By \cite[3.15(1)]{GP2} the canonical homomorhisms
\[
H_0(\ZF(\Delta^{\bullet}\times \Gm\times\Gm,Y)) \to
H_0(\ZF(\Delta^{\bullet}\times \mathbb G_{m,k(u)},Y))
\]
and
\[
H_0(\ZF(\Delta^{\bullet}\times \mathbb G_{m,k(u)},Y)) \to
H_0(\ZF(\Delta^{\bullet}_{\Spec k(t,u)},Y))
\]
are injective, hence the lemma.
\end{proof}

Let $Y$ be a $k$-smooth variety and $F/k$ be a field extension.
There is a map of pointed sets
   $$adj:\Fr_n(\Spec(F),Y)\to \Fr^F_n(\Spec(F),Y_F)$$
given by the assignment $(Z,W,\phi,g)\mapsto (Z,W,\phi^F,g^F)$.
Here for a $k$-morphism $g: W\to Y$ we write
$g^F$ to denote the $F$-morphism
$(g, pr_{\text{Spec}(F)}): W\to Y_F$ and
$pr_{\text{Spec}(F)}: W\to \text{Spec}(F)$
is the structure morphism.
In particular, for $Y=\A^n$ and $\phi: W\to \A^n_k$
we write $\phi^F$ for $(\phi,pr_{\text{Spec}(F)}): W\to \A^n_F$.
It is easy to see that the map $adj$ is a bijection. Moreover,
it induces bijections
   $$adj: \ZF_n(\Spec(F),Y)\to \ZF^F_n(\Spec(F),Y_F) \ \text{and} \ \ZF(\Spec(F),Y)\to \ZF^F(\Spec(F),Y_F).$$

\begin{lemma} \label{l:switch}
Let $F/k$ be a field extension, choose $x,y\in F^{\times}$ such that $x\neq
y^{\pm 1}$ and let $u_1,u_2$ be coordinates on $\Gm\times\Gm$.
Consider morphisms $f,g\colon\Spec F\to \Gm\times \Gm$ given by
$u_1\mapsto x, u_2\mapsto y$ and $u_1\mapsto y, u_2\mapsto x$
respectively. Then for $p=(\id-e_1)\boxtimes (\id-e_1)$ we have
$$p\circ (f\boxtimes \Sigma) \sim p\circ (g\boxtimes (-\varepsilon))$$
in $\ZF(\Spec F,\Gm\times \Gm)$.
\end{lemma}

\begin{proof}
The above adjunction isomorphism
\[
adj: \ZF(\Spec F,\Gm\times \Gm)\cong \ZF^F(\Spec F,\Gm_{,F}\times
\Gm_{,F})
\]
implies it is sufficient to verify the case $F=k$. So we have
morphisms $f,g\colon \pt \to\Gm\times \Gm$, $\pt\mapsto (x,y)$ and $\pt
\mapsto (y,x)$ respectively. Taking suspensions, we obtain framed
correspondences
\[
(\{0\}, \A^1, t, c_{(x,y)}),\,(\{0\}, \A^1, t, c_{(y,x)})\in
\Fr_1(\pt,\Gm\times\Gm),
\]
where $c_{(x,y)}$ and $c_{(y,x)}$ are morphisms on $\A^1$ sending it
to the points $(x,y)$ and $(y,x)$ respectively.

Consider $h(s,t)=\frac{1}{x-y}(t^2-(s(x+y)+(1-s)(xy+1))t+xy) \in
k[s,t,t^{-1}]=k[\A^1\times \Gm]$ and a framed correspondence
\begin{equation}\label{e:switching}
H_s:=(Z(h), \A^1\times \Gm, h(s,t), (t,xyt^{-1})\circ pr_{\bb G_m})\in
\Fr_1(\A^1,\Gm\times\Gm).
\end{equation}
We have $h(0,t)=\frac{1}{x-y}(t-xy)(t-1)$ and
$h(1,t)=\frac{1}{x-y}(t-x)(t-y)$. Using the additivity property for
supports in $\ZF_1(\pt,\Gm\times\Gm)$ (see Definition~\ref{stab})
and Lemma~\ref{l:canonicalpolynomial} we will check below that
\begin{equation}\label{r:switching}
(\{0\}, \A^1, t, c_{(x,y)})+(\{0\}, \A^1, -t, c_{(y,x)}) \sim
(\{0\}, \A^1, \frac{1-xy}{x-y}t, c_{(1,xy)})+ (\{0\}, \A^1,
\frac{xy-1}{x-y}t, c_{(xy,1)})
\end{equation}
in $\ZF_1(\pt,\Gm\times\Gm)$. The composition with the idempotent
$p$ annihilates all extra summands and proves the lemma.

In order to prove the relation (\ref{r:switching}), consider the
framed correspondence (\ref{e:switching}) in
$\ZF_1(\A^1,\Gm\times\Gm)$. Observe that in
$\ZF_1(\pt,\Gm\times\Gm)$
\begin{multline*}
H_1=(Z(h(t,1),\Gm,h(1,t),(t,xyt^{-1}))= \\
=(\{x\},\Gm-\{y\},\frac{1}{x-y}(t-x)(t-y),(t,xyt^{-1}))+(\{y\},\Gm-\{x\},\frac{1}{x-y}(t-x)(t-y),(t,xyt^{-1})).
\end{multline*}
By Lemma~\ref{l:canonicalpolynomial} one has in
$\ZF_1(\pt,\Gm\times\Gm)$
$$(\{x\},\Gm-\{y\},\frac{1}{x-y}(t-x)(t-y),(t,xyt^{-1}))\sim (\{0\}, \A^1, \frac{x-y}{x-y}t, c_{(x,y)})=(\{0\}, \A^1, t, c_{(x,y)}),$$
$$(\{y\},\Gm-\{x\},\frac{1}{x-y}(t-x)(t-y),(t,xyt^{-1}))\sim (\{0\}, \A^1, \frac{y-x}{x-y}t, c_{(x,y)})=(\{0\}, \A^1, -t, c_{(x,y)}).$$
Thus $H_1\sim (\{0\}, \A^1, t, c_{(x,y)})+(\{0\}, \A^1, -t,
c_{(y,x)})$ in $\ZF_1(\pt,\Gm\times\Gm)$. Similar computations show
that $H_0\sim (\{0\}, \A^1, \frac{1-xy}{x-y}t, c_{(1,xy)})+ (\{0\},
\A^1, \frac{xy-1}{x-y}t, c_{(xy,1)})$ in $\ZF_1(\pt,\Gm\times\Gm)$.
The equality (\ref{r:switching}) is proved.
Since the right hand side of the equality
(\ref{r:switching}) is annihilated by the idempotent $p$, our lemma follows.
\end{proof}

\begin{proposition} \label{p:switch}
Let $\tau\colon \Gm\times \Gm \to \Gm\times \Gm$ be the permutation
of coordinates morphism. Denote $p=(\id-e_1)\boxtimes (\id-e_1)$.
Then $p \circ (\id_{\Gm\times \Gm}\boxtimes (-\epsilon)) \sim p\circ (\tau\boxtimes \Sigma)$ in
$\ZF(\Gm\times \Gm,\Gm\times \Gm)$.
\end{proposition}

\begin{proof}
Note that $\tau\boxtimes (-\varepsilon)=(\id_{\Gm\times \Gm}\boxtimes (-\varepsilon))\circ \tau$. Hence
$$p\circ (\tau\boxtimes (-\varepsilon))\circ \tau=p \circ (\id_{\Gm\times \Gm}\boxtimes (-\epsilon)) \in \ZF_1(\Gm\times \Gm, \Gm\times \Gm)).$$
Similarly, $p\circ (\id_{\Gm\times \Gm}\boxtimes \Sigma)\circ\tau=p \circ (\tau\boxtimes \Sigma)$ in $\ZF_1(\Gm\times \Gm, \Gm\times \Gm))$.
It remains to check that
   \begin{equation}\label{sw_2}
    p\circ (\id_{\Gm\times \Gm}\boxtimes \Sigma) = p\circ (\tau\boxtimes (-\varepsilon))
   \end{equation}
in $H_0(\ZF(\Delta^{\bullet}\times \Gm\times \Gm,\Gm\times \Gm))$.

Let $u_1$ and $u_2$ be coordinate functions on $\Gm\times \Gm$. Let
$f\colon \Spec k(u_1,u_2)\to \Spec k[u_1,u_2]$ be the canonical embedding and
$g\colon \Spec k(u_1,u_2)\to \Spec k[u_1,u_2]$ be given by $g^*(u_1)=u_2, g^*(u_2)=u_1$.
By Proposition \ref{l:switch} we know that
$p\circ (f\boxtimes \Sigma) \sim p\circ (g\boxtimes (-\varepsilon))$ in
$H_0(\ZF(\Delta^{\bullet}_{k(u_1,u_2)},\Gm\times \Gm))$.
Proposition~\ref{p:injectivity} yields the desired equality~\eqref{sw_2} in
$H_0(\ZF(\Delta^{\bullet}\times \Gm\times \Gm,\Gm\times \Gm))$.
\end{proof}


Recall that $\Sigma=(\{0\},\A^1,t)\in\ZF_1(\pt,\pt)$. For every $k>0$
we write $\Sigma^k$ to denote $\Sigma\boxtimes\bl
k\cdots\boxtimes\Sigma\in\ZF_k(\pt,\pt)$.

Let $\tau\colon \Gm\times \Gm\to \Gm\times \Gm$
be the permutation of coordinates morphism.
For each even integer $m\geq 0$ and each integer $n\geq 0$ consider two presheaf morphisms
$$(-\boxtimes \Sigma^{2n}): \ZF_m(-\times X\times \Gm\times \Gm, Y \times \Gm\times \Gm) \to \ZF_{m+2n}(-\times X\times \Gm\times \Gm, Y \times \Gm\times \Gm),$$
$$(-\boxtimes \Sigma^{2n})\circ sw:\ZF_m(-\times X\times \Gm\times \Gm, Y \times \Gm\times \Gm) \to \ZF_{m+2n}(-\times X\times \Gm\times \Gm, Y \times \Gm\times \Gm),$$
where $sw(a)=(\id_Y\times \tau)\circ a\circ (\id_X\times \tau)$.

\begin{lemma} \label{l:switch2}
Let $X,Y$ be $k$-smooth schemes. Given an even integer $m\geq 0$, there exists a large enough
$n$ and a homotopy
$$H: \ZF_m(-\times X\times \Gm\times \Gm, Y \times \Gm\times \Gm) \to 
\ZF_{m+2n}(-\times X\times \Gm\times \Gm\times \A^1, Y \times \Gm\times \Gm)$$
such that for any $a\in \ZF_m(W\times X\wedge (\Gm,1)\wedge
(\Gm,1),Y\wedge (\Gm,1)\wedge (\Gm,1))$ one has
\[
a\boxtimes \Sigma^{2n}=H_0(a) \ \ \text{and} \ \
H_1(a)=\Sigma^{2n}([(\id_Y\times \tau)\circ a\circ (\id_X\times
\tau)]).
\]
And both $H_0(a)$ and $H_1(a)$ are in $\ZF_{m+2n}(W\times
X\wedge (\Gm,1)\wedge (\Gm,1), Y\wedge (\Gm,1)\wedge
(\Gm,1))$.
\end{lemma}

\begin{proof}
It follows from Proposition~\ref{p:switch} that there exists
a large enough integer $n$ and a homotopy
$\Psi\in \ZF_n (\Gm\times \Gm \times \A^1,\Gm\times \Gm)$
such that
$i_0^*(\Psi)=p\circ(-\varepsilon \boxtimes
(\Sigma^{n-1}\id_{\Gm\times \Gm }))$
and
$i_1^*(\Psi)=p\circ \Sigma^{n}\tau$,
where
$p=(\id-e_1)\boxtimes (\id-e_1)$.

Given any element $a\in \ZF_m(W\times X\times \Gm\times \Gm, Y\times
\Gm\times \Gm)$, set
\begin{multline*}
H'(a)=(id_Y\times \Psi)\circ (a\times id_{\A^1})\circ (id_{W\times X}\times \Psi\times id_{\A^1})\circ (id_{W\times X\times \Gm \times \Gm}\times \Delta) \in \\
\in \ZF_{m+2n}(W\times X\times \Gm\times \Gm\times \A^1, Y\times
\Gm\times \Gm),
\end{multline*}
where $\Delta: \A^1 \to \A^1\times\A^1$ is the diagonal morphism.
Then for any element $a\in \ZF_m(W\times X\wedge (\Gm,1)\wedge
(\Gm,1),Y\wedge (\Gm,1)\wedge (\Gm,1))$ one has
$$H'(a)_0=[id_Y \times \Sigma^{n-1}(-\epsilon)]\circ a\circ[id_X \times \Sigma^{n-1}(-\epsilon)] \ \ \text{and}
\ \ H'(a)_1=[\id_Y\times \Sigma^n(\tau)]\circ a\circ[\id_X\times
\Sigma^n(\tau)].$$ It is easy to see that there are matrices $A,B\in
SL_{m+2n}(k)$ such that for any element $a$ in $\ZF_m(W\times
X\wedge (\Gm,1)\wedge (\Gm,1),Y\wedge (\Gm,1)\wedge (\Gm,1))$ one
has
\begin{multline*}
\varphi_A([id_Y \times \Sigma^{n-1}(-\epsilon)]\circ  a \circ [id_X \times \Sigma^{n-1}(-\epsilon)])=a\boxtimes \Sigma^{2n}=\Sigma^{2n}(a), \\
\varphi_B([\id_Y\times \Sigma^n(\tau)]\circ a\circ[\id_X\times
\Sigma^n(\tau)])= ([\id_Y\times \tau]\circ a\circ[\id_X\times
\tau))\boxtimes \Sigma^{2n}= \Sigma^{2n}([\id_Y\times \tau]\circ
a\circ[\id_X\times \tau]).
\end{multline*}
Choose matrices $A_s,B_s \in SL_{m+2n}(k[s])$ such that $A_0=id$,
$A_1=A$, $B_0=id$, $B_1=B$. Then for the matrix $C_s=B_s\circ
A_{1-s} \in SL_{m+2n}(k[s])$ one has $C_0=A$, $C_1=B$. Set
$H=\varphi_{C_s} \circ H'$. Then for the chosen element $a \in
\ZF_m(W\times X\wedge (\Gm,1)\wedge (\Gm,1), Y\wedge (\Gm,1)\wedge
(\Gm,1))$, one has
$$H_0(a)=\varphi_{A}(H'(a)_0)=\Sigma^{2n}(a) \ \ \text{and}
\ \ H_1(a)=\varphi_{B}(H'(a)_1)=\Sigma^{2n}([\id_Y\times \tau]\circ
a\circ[\id_X\times \tau)),$$ as was to be proved.
\end{proof}

\section{The inverse morphism}

The main aim of this section is to define for any integers $n,m\geq
0$ a subpresheaf $\ZF^{(n)}_m(-\times \Gm,Y\times\Gm)$ of the
presheaf $\ZF_m(-\times \Gm,Y\times \Gm)$ and define a morphism of
abelian presheaves
   $$\rho_n: \ZF^{(n)}_m(-\times \Gm,Y\times \Gm) \to \ZF_m(-,Y).$$
We also prove certain properties of morphisms $\rho_n$ and of
presheaves $\ZF^{(n)}_m(-\times \Gm,Y\times \Gm)$ which are used in
the proof of the Linear Cancelation Theorem (Theorem~C).

We begin with some general remarks. Let $X$ and $Y$ be $k$-smooth
schemes. Consider a framed correspondence
\[
a=(Z,U,(\varphi_1,\varphi_2,\hdots,\varphi_m),g)\in \Fr_m(X\times
\Gm,Y\times \Gm).
\]
Let $(U, p:U\to \A^m\times (X\times\Gm),s: Z\to U)$ be the \'{e}tale
neighborhood of $Z$ in $\A^m\times (X\times\Gm)$ from the definition
of the framed correspondence $a$. Let $t$ be the invertible function
on $X\times\Gm$ corresponding to the projection on $\Gm$ and $u$ be
invertible function on $Y\times\Gm$ corresponding to the projection
on $\Gm$. Let $f_2=g^*(u)$ and $f_1=p^*_{X\times\Gm}(t)$ be two
invertible functions on $U$, where
$p_{X\times\Gm}=pr_{X\times\Gm}\circ p: U \to X\times\Gm$. Set
$g=(g_1,g_2)$, where $g_1=pr_Y\circ g$ and $g_2=pr_{\Gm}\circ g$.

Since $Z$ is finite over $X\times \Gm$, the $\mathcal
O_{X\times\Gm\times Y\times\Gm}$-sheaf $P_a:=\mathcal
O_U/(\varphi_1,\varphi_2,\hdots,\varphi_m)$ is finite over
$X\times\Gm$. Since the sheaf $P_a$ is finite over $X\times
\Gm$,
it is flat over $X\times\Gm$
by \cite[Lemma 7.3]{OP}.

Let $Z^+_n$ be the closed subset of $Z$ defined by the equation
$(f^{n+1}_1-1)|_Z=0$. Let $Z^-_n$ be the closed subset of $Z$
defined by the equation $(f^{n+1}_1-f_2)|_Z=0$. Note that $Z^+_n$ is
finite over $X$ if and only if $\mathcal
O_U/(f^{n+1}_1-1,\varphi_1,\varphi_2,\hdots,\varphi_m)$ is finite
over $X$. By \cite[4.1]{S} the latter $\mathcal O_X$-module is
always finite and even flat. Note that $Z^-_n$ is finite over $X$ if
and only if $\mathcal
O_U/(f^{n+1}_1-f_2,\varphi_1,\varphi_2,\hdots,\varphi_m)$ is finite
over $X$. 
We have mentioned above that the $\mathcal O_{X\times \Gm}$-module
$P_a=\mathcal O_U/(\varphi_1,\varphi_2,\hdots,\varphi_m)$ is finite
and flat over $\mathcal O_{X\times \Gm}$. Thus by \cite[4.1.b]{S} 
there exists an integer $N$ such that for any $n\geq N$
the $\mathcal O_X$-module
$\mathcal O_U/(f^{n+1}_1-f_2,\varphi_1,\varphi_2,\hdots,\varphi_m)$
is finite and even flat over $X$. 
In particular, $Z^-_n$ is finite over $X$ for any $n\geq N$.

The following definition is inspired by \cite[Section 4]{S}.
\begin{definition}{\rm
Let $X$ and $Y$ be $k$-smooth schemes. Consider a framed
correspondence $a=(Z,U,(\varphi_1,\varphi_2,\hdots,\varphi_m),g)\in
\Fr_m(X\times \Gm,Y\times \Gm)$. Set
   $$\rho^+_{n,fr}(a):=(Z^+_n,U,(f^{n+1}_1-1,\varphi_1,\varphi_2,\hdots,\varphi_m),g_1)$$
and
   $$\rho^-_{n,fr}(a):=(Z^-_n,U,(f^{n+1}_1-f_2,\varphi_1,\varphi_2,\hdots,\varphi_m),g_1).$$
As we have mentioned above, $Z^+_n$ is finite over $X$ for all
$n\geq 0$, hence $\rho^+_{n,fr}(a)\in \ZF_{m+1}(X,Y)$. We say that
$\rho^-_{n,fr}(a)$ is {\it defined\/} if $Z^-_n$ is finite over $X$,
which is equivalent to saying that the $\mathcal O_X$-module
$P_a/(f^{n+1}_1-f_2)P_a$ is finite and flat over $X$. If
$\rho^-_{n,fr}(a)$ is defined, then we set
$$\rho_{n,fr}(a)=\rho^+_{n,fr}(a)-\rho^-_{n,fr}(a) \in \ZF_{m+1}(X,Y)$$
and say that $\rho_{n,fr}(a)$ is {\it defined}.
}\end{definition}

Given integers $m,n\geq 0$, denote by $\Fr^{(n)}_m(X\times
\Gm,Y\times \Gm)$ the subset of those framed correspondences $a\in
\Fr_m(X\times \Gm,Y\times \Gm)$ for which the $\mathcal O_X$-module
$P_a/(f^{n+1}_1 -f_2)P_a$ is finite over $X$ (that is
$\rho_{n,fr}(a)$ is defined). It follows from~\cite[4.4]{S} that the
assignment $X'\mapsto \Fr^{(n)}_m(X'\times \Gm,Y\times \Gm)$ is a
subpresheaf of $\Fr_m(-\times \Gm,Y\times \Gm)$.

\begin{definition}\label{d:presheaf}{\rm
Define a presheaf of abelian groups $\ZF^{(n)}_m(-\times \Gm,Y\times
\Gm)$ as follows. Its sections on $X$ is the abelian group $\mathbb
Z[\Fr^{(n)}_m(X\times \Gm,Y\times \Gm)]$ modulo a subgroup generated
by all elements of the form
$$(Z_1\sqcup Z_2,U_1\sqcup U_2,\varphi_1\sqcup \varphi_2,g_1\sqcup g_2)-(Z_1,U_1,\varphi_1,g_1)-(Z_2,U_2,\varphi_2,g_2).$$
It is straightforward to check that $\ZF^{(n)}_m(X\times \Gm,Y\times
\Gm)$ is a free abelian group with a free basis consisting of the
elements of the form $a=(Z,U,\varphi,g)$, where $Z$ is connected and
the $\mathcal O_X$-module $P_a/(f^{n+1}_1 -f_2)P_a$ is finite and
flat over $X$. Moreover, the group $\ZF^{(n)}_m(X\times \Gm,Y\times
\Gm)$ is a subgroup of the group $\ZF_m(X\times \Gm,Y\times \Gm)$,
and $\ZF^{(n)}_m(-\times \Gm,Y\times \Gm)$ is a subpresheaf of the
presheaf $\ZF_m(-\times \Gm,Y\times \Gm)$.

}\end{definition}

It follows from~\cite[4.4]{S} that for any morphism $f: X' \to X$ of
smooth varieties the following diagram is commutative
$$\xymatrix{\Fr^{(n)}_m(X\times \Gm,Y\times \Gm) \ar^{(f\times id)^*}[r]\ar^{\rho_{n,fr}}[d]& \Fr^{(n)}_m(X'\times \Gm,Y\times \Gm) \ar^{\rho_{n,fr}}[d]\\
              \ZF_{m+1}(X,Y)  \ar^{f^*}[r]& \ZF_{m+1}(X',Y).}$$
We see that $\rho_{n,fr}:\Fr^{(n)}_m(-\times \Gm,Y\times \Gm) \to
\ZF_{m+1}(-,Y)$ is a morphism of pointed presheaves. We can extend
it to get a morphism of presheaves of abelian groups $\mathbb
Z[\Fr^{(n)}_m(-\times \Gm,Y\times \Gm)] \to \ZF_{m+1}(-,Y)$. This
morphism annihilates the elements of the form
$$(Z_1\sqcup Z_2,U_1\sqcup U_2,\varphi_1\sqcup \varphi_2,g_1\sqcup g_2)-(Z_1,U_1,\varphi_1,g_1)-(Z_2,U_2,\varphi_2,g_2).$$

\begin{definition}\label{d:presheaf_and_morphism}{\rm
The above arguments show that the presheaf morphism $\rho_{n,fr}$
induces a unique presheaf of abelian groups morphism
   $$\rho_n: \ZF^{(n)}_m(-\times \Gm,Y\times \Gm) \to \ZF_{m+1}(-,Y)$$
such that for any
$a \in Fr^{(n)}_m(X\times \Gm,Y\times \Gm)$
one has
$\rho_n(a)=\rho_{n,fr}(a)$. We also call $\rho_n$ the {\it
inverse morphism}.

}\end{definition}

\begin{lemma}\label{l:colimit}
The following relations are true:
$$\Fr_m(-\times \Gm,Y\times \Gm)=\colim_n \Fr^{(n)}_m(-\times \Gm,Y\times \Gm),$$
$$\ZF_m(-\times \Gm,Y\times \Gm)=\colim_n \ZF^{(n)}_m(-\times \Gm,Y\times \Gm).$$
\end{lemma}

This lemma follows from the following

\begin{proposition}(\cite[4.1]{S}) For any framed correspondence $a \in \Fr_m(X\times \Gm,Y\times \Gm)$ one
has:

\begin{itemize}
\item[(a)] for any $n \geq 0$, the sheaf $P_a/(f^{n+1}_1 -1)P_a$ is
finite and flat over X;
\item[(b)] there exists an integer $N$ such
that, for any $n\geq N$, the sheaf $P_a/(f^{n+1}_1 -f_2)P_a$ is
finite and flat over $X$.
\end{itemize}
\end{proposition}

We shall need the following obvious property of $\rho_n$.

\begin{lemma}\label{rho_and_suspension}
For any integers $m,n,r \geq 0$, the following diagram commutes
$$\xymatrix{\ZF^{(n)}_m(-\times \Gm,Y\times \Gm) \ar^{\Sigma^r}[r]\ar^{\rho_n}[d]& \ZF^{(n)}_{m+r}(-\times \Gm,Y\times \Gm) \ar^{\rho_n}[d]\\
              \ZF_{m+1}(-,Y)  \ar^{\Sigma^r}[r]& \ZF_{m+1+r}(-,Y).}$$
\end{lemma}

\begin{lemma}\label{l:a_boxtimes}
Let $X$ and $Y$ be $k$-smooth schemes. Then for any integers $m$ and
$n$ and any $a\in \ZF_m(X,Y)$, one has $a\boxtimes (\id-e_1) \in
\ZF^{(n)}_m(X\times\Gm,Y\times\Gm)$. In particular, for any integers
$m$ and $n$ there is a well defined composite morphism
   $$\rho_n\circ (- \boxtimes (\id-e_1)): \ZF_m(-\times X,Y)\to \ZF^{(n)}_m(-\times X\times\Gm,Y\times\Gm) \to \ZF_{m+1}(-\times X,Y).$$
Moreover, for an element $a \in \ZF_m(W\times X,Y)$ of the form
$(Z,U,(\varphi_1,\varphi_2,\dots,\varphi_{m}),g)$ one has
\begin{multline*}
\rho_n(a\boxtimes (\id-e_1))=-(Z\times Z(t^{n+1}-t),U\times\Gm,(t^{n+1}-t,\varphi_1,\varphi_2,\dots,\varphi_{m}),g) +\\
+ (Z\times
Z(t^{n+1}-1),U\times\Gm,(t^{n+1}-1,\varphi_1,\varphi_2,\dots,\varphi_{m}),g)
\in \ZF_{m+1}(W\times X,Y).
\end{multline*}
\end{lemma}

\begin{proof}
Let $a \in \ZF_m(W\times X,Y)$ be the image of
$(Z,U,(\varphi_1,\varphi_2,\dots,\varphi_{m}),g)\in \Fr_{m}(W\times
X,Y).$ Then
\begin{multline*}
a\boxtimes (\id-e_1)=(Z\times \Gm,U\times\Gm,(\varphi_1,\varphi_2,\dots,\varphi_{m}),(g,t)) -\\
- (Z\times
\Gm,U\times\Gm,(\varphi_1,\varphi_2,\dots,\varphi_{m}),(g,e_1)) \in
\ZF_{m}(W\times X\times \Gm,Y\times \Gm),
\end{multline*}
where $t$ is the coordinate function on $\Gm$. Clearly,
$Z^+_n=Z\times Z(t^{n+1}-1)\subset Z\times\Gm$ and  $Z^-_n=Z\times
Z(t^{n+1}-t)\subset Z\times\Gm$. Both sets are finite over $X$.
Hence $a\boxtimes (\id-e_1) \in \ZF^{(n)}_m(X\times\Gm,Y\times\Gm)$
in this case. Any element of $\ZF_m(W\times X,Y)$ is a linear
combination of elements of the form
$(Z,U,(\varphi_1,\varphi_2,\dots,\varphi_{m}),g)$. This proves the
first assertion of the lemma.

Computing $\rho_n(a\boxtimes (\id-e_1))$ for
$a=(Z,U,(\varphi_1,\varphi_2,\dots,\varphi_{m}),g)$ we obtain
\begin{multline*}
\rho_n(a\boxtimes (\id-e_1))=-(Z\times Z(t^{n+1}-t),U\times\Gm,(t^{n+1}-t,\varphi_1,\varphi_2,\dots,\varphi_{m}),g) +\\
+ (Z\times
Z(t^{n+1}-1),U\times\Gm,(t^{n+1}-1,\varphi_1,\varphi_2,\dots,\varphi_{m}),g)
\in \ZF_{m+1}(W\times X,Y),
\end{multline*}
as was to be shown.
\end{proof}

\begin{lemma} \label{l:composition}
Let $X$ and $Y$ be $k$-smooth schemes.
Then for every even integer $m$ and any $n$ one has
\[
\rho_n\circ (- \boxtimes (\id-e_1))\sim (-\boxtimes \epsilon):
\ZF_m(-\times X,Y)\rightrightarrows \ZF_{m+1}(-\times X,Y),
\]
where $\epsilon=(\{0\},\A^1, -t,c') \in \ZF_1(\pt,\pt)$.
\end{lemma}

\begin{proof}
Set $\eta_n=\rho_n \circ (- \boxtimes (\id-e_1))$. Take the matrix
\[
A=\begin{pmatrix}
0 & 1 & 0 & \hdots & 0 & 0\\
0 & 0 & 1 & \hdots & 0 & 0\\
\vdots & \vdots  & \vdots  & \ddots & \vdots  & \vdots \\
0 & 0 & 0 & \hdots & 0 & 1\\
1 & 0 & 0 & \hdots & 0 & 0\\
\end{pmatrix}\in SL_{m+1}(k)
\]
and let $A_s\in SL_{m+1}(k[s])$ be such that $A_0=id$, $A_1=A$. Let
$H_{A_s}$ be the $\A^1$-homotopy from Definition~\ref{d:homotopy_1}
between the identity and $\varphi_A$. By
Definition~\ref{d:homotopy_1} one has
   $$\eta_n=\rho_n \circ (- \boxtimes (\id-e_1)) \frac{H_{A_s}\circ \eta_n}{} \varphi_A \circ \rho_n \circ (- \boxtimes (\id-e_1))= \varphi_A\circ \eta_n.$$
Set $H'=H_{A_s}\circ \eta_n$. By Lemma \ref{l:H1_and_H_2} it remains
to find an $H''$ such that $\varphi_A \circ \eta_n\frac{H''}{}
(-\boxtimes \epsilon)$ and set $H=H'+H''-H_{\varphi_A\circ \eta_n}$.
In this case by Lemma \ref{l:H1_and_H_2} one gets $\rho_n \circ (-
\boxtimes (\id-e_1))=\eta_n\frac{H}{} (-\boxtimes \epsilon)$.

To construct $H''$, note that by the last statement of Lemma
\ref{l:a_boxtimes} one has
$$\varphi_A \circ \eta_n = -\boxtimes [(Z(t^{n+1}-1),\Gm, t^{n+1}-1,c)-(Z(t^{n+1}-t),\Gm, t^{n+1}-t,c)]$$
and $(-\boxtimes \epsilon)=-\boxtimes (\{0\},\A^1,-t,c')$,
where $c\colon \Gm\to \pt$ is the canonical projection. By Lemma
\ref{l:box_b1_and_b2} one can take $H''$ to be an $\A^1$-homotopy of
the form $H''=(-\boxtimes h'')$, where $h'' \in \ZF_1(\A^1,pt)$ is
such that
$$(Z(t^{n+1}-1),\Gm, t^{n+1}-1,c)- (Z(t^{n+1}-t),\Gm,
t^{n+1}-t,c)=h''_0$$
and
$$h''_1=(\{0\},\A^1,-t,c') \in
\ZF_1(\pt,\pt),$$ where $c'\colon \A^1\to \pt$ is the canonical
projection. Now let us find the desired element $h''$. Since
$t^{n+1}-1$ does not vanish at $t=0$, we can extend the neighborhood
from $\Gm$ to $\A^1$ to get an equality,
\[
(Z(t^{n+1}-1),\Gm, t^{n+1}-1,c)=(Z(t^{n+1}-1),\A^1, t^{n+1}-1),c')
\in \ZF_1(pt,pt).
\]
By Lemma~\ref{l:polynomialcycle} there is $h''' \in \ZF_1(\A^1,pt)$
such that
\[
(Z(t^{n+1}-1),\A^1, t^{n+1}-1,c')=h'''_0 \ \ \text{and} \ \
h'''_1=(Z(t^{n+1}-t),\A^1, t^{n+1}-t,c') \in \ZF_1(pt,pt),
\]
because polynomials $t^{n+1}-t$ and $t^{n+1}-1$ have the same degree
and the same leading coefficient. Using the additivity property for
supports in $\ZF_1(pt,pt)$ and the second statement of
Lemma~\ref{l:canonicalpolynomial}, we can find an element $h^{iv}\in
\ZF_1(\A^1,pt)$ such that
\[
(Z(t^{n+1}-t),\Gm, t^{n+1}-t,c)=h^{iv}_0  \ \ \text{and} \ \
h^{iv}_1=(Z(t^{n+1}-t),\A^1, t^{n+1}-t,c')-(\{0\},\A^1,-t,c') \in
\ZF_1(pt,pt)
\]
Set $h'':=h'''- h^{iv} \in \ZF_1(\A^1,pt)$.
Then $h''$ is the desired element.

Set $H''=(-\boxtimes h'')$ and $H=H'+H''-H_{\varphi_A\circ \eta_n}$.
Then $H$ is the desired $\A^1$-homotopy. That is
$$\rho_n \circ (- \boxtimes (\id-e_1)) \frac{H}{} (-\boxtimes \epsilon)$$
and our statement follows.
\end{proof}

\section{Theorem C}

The main purpose of this section is to prove Theorem~C. We sometimes
identify simplicial abelian groups with chain complexes concentrated
in non-negative degrees by using the Dold--Kan correspondence.

\begin{lemma}\label{l:Moore_1}
Let $X$ and $Y$ be $k$-smooth schemes and $m,r,N\geq 0$ be integers.
Then for any Moore cycle $a \in \ZF_m(\Delta^r \times X,Y)$ of the
simplicial abelian group $\ZF_m(\Delta^\bullet\times X,Y)$, one has
$a\boxtimes (id-e_1) \in \ZF^{(N)}_{m}(\Delta^r \times X\times
\Gm,Y\times \Gm)$. Moreover, $\rho_N(a\boxtimes (id-e_1))$ is a
Moore cycle. The homology classes of Moore cycles
   $$a\boxtimes \epsilon \ \ \text{and} \ \ \rho_N(a\boxtimes (id-e_1))$$
coincide in $\ZF_{m+1}(\Delta^{\bullet} \times X,Y)$.
\end{lemma}

\begin{proof}
The element $a\boxtimes (id-e_1)$ is in $\ZF^{(N)}_{m}(\Delta^r
\times X\times \Gm,Y\times \Gm)$ by Lemma \ref{l:a_boxtimes}. Since
$\ZF^{(N)}_m(-\times \Gm,Y\times \Gm ))$ is a presheaf, then
$\partial_i(a\boxtimes (id-e_1))\in \ZF^{(N)}_{m}(\Delta^r \times
X\times \Gm,Y\times \Gm)$. Since the morphism $\rho_N$ is a morphism
of presheaves, then
$$\partial_i(\rho_N(a\boxtimes (id-e_1)))=\rho_N(\partial_i(a\boxtimes (id-e_1))=\rho_N(\partial_i(a)\boxtimes (id-e_1))=0.$$
This proves the first assertion of the lemma.

By Lemma \ref{l:composition} the morphism
$$a'\mapsto \rho_N(a'\boxtimes (\id_{\Gm}-e_1)): \ZF_m(- \times X,Y) \to \ZF^{(N)}_m(- \times X,Y) \to \ZF_{m+1}(-\times X,Y)$$
is $\A^1$-homotopic to the morphism $a'\mapsto a'\boxtimes
\epsilon$. Thus the corresponding morphisms of the simplicial
abelian groups $\ZF_m(\Delta^\bullet \times X,Y) \rightrightarrows
\ZF_{m+1}(\Delta^\bullet \times X,Y)$ induce the same morphisms on
homology. Hence the homology class of the Moore cycle
$\rho_N(a\boxtimes (\id_{\Gm}-e_1))$ coincides with the homology
class of the Moore cycle $a\boxtimes \epsilon$.
\end{proof}

\begin{lemma}\label{espsilon_and_sigma}
One has $\epsilon\boxtimes \epsilon \sim \Sigma^2$ in
$\ZF_2(\pt,\pt)$. Moreover, for any integer $r\geq 0$ one has
$\epsilon\boxtimes \epsilon\boxtimes \Sigma^{r}\sim \Sigma^{2+r}$ in
$\ZF_{2+r}(\pt,\pt)$.
\end{lemma}
\begin{proof}
Let $c: \A^1\times \A^2 \to \pt$ be the structure morphism. Take the matrix
\[
A=\begin{pmatrix}
0 & -1 \\
-1 & 0 \\
\end{pmatrix}\in SL_{2}(k).
\]
There is an $A_s\in SL_{2}(k[s])$ such that $A_0=id$, $A_1=A$. Take
$$h_s=(\A^1\times 0,\A^1\times \A^2, A_s\circ (t_1,t_2),c) \in \ZF_2(\A^1,\pt).$$
Clearly, $h_0=\Sigma^2$ and $h_1=\epsilon\boxtimes \epsilon$.
The first assertion is proved.
To prove the second one take the element
$h_s\boxtimes \Sigma^r \in \ZF_{2+r}(\A^1,\pt)$. Then
$h_0\boxtimes \Sigma^r=\Sigma^{2+r}$
and
$h_1\boxtimes \Sigma^r=\epsilon\boxtimes \epsilon\boxtimes \Sigma^r$.
\end{proof}

\begin{corollary}\label{box_epsilon_2_and_Sigma_2}
Let $X$ and $Y$ be $k$-smooth schemes and $m\geq 0$ be an integer.
Then,
$$(-\boxtimes \epsilon^2)\sim (-\boxtimes \Sigma^2): \ZF_m(-\times X,Y)\rightrightarrows \ZF_{m+2}(-\times X,Y)$$
and
$$(-\boxtimes \epsilon^2\boxtimes \Sigma^{r})\sim (-\boxtimes \Sigma^{2+r}): \ZF_m(-\times X,Y)\rightrightarrows \ZF_{m+2+r}(-\times X,Y).$$
Therefore the first pair of maps produces the same maps on homology
$$H_*(\ZF_m(\Delta^{\bullet}\times X,Y))\rightrightarrows H_*(\ZF_{m+2}(\Delta^{\bullet}\times X,Y)).$$
Similarly, the second pair of maps gives the same maps on homology
$$H_*(\ZF_m(\Delta^{\bullet}\times X,Y))\rightrightarrows H_*(\ZF_{m+2+r}(\Delta^{\bullet}\times X,Y)).$$
\end{corollary}

\begin{lemma}\label{l:injectivity}
Let $X$ and $Y$ be $k$-smooth schemes and $m\geq 0$ be an integer.
Then for any integer $r\geq 0$ one has
\begin{multline*}
\text{Ker}[-\boxtimes (\id_\Gm - e_1):H_r(\ZF_m(\Delta^\bullet \times X,Y)) \to H_r(\ZF_m((\Delta^\bullet \times X)\wedge (\Gm,1),Y\wedge (\Gm,1)))]\subseteq \\
\subseteq \text{Ker}[(-\boxtimes \Sigma^2):H_r(\ZF_m(\Delta^\bullet
\times X,Y))\to H_r(\ZF_{m+2}(\Delta^\bullet \times X,Y))].
\end{multline*}
\end{lemma}

\begin{proof}
Consider the associated Moore complexes. Assume that
\[
a\in \ZF_m(\Delta^r\times X,Y)
\]
is a Moore cycle for which $a\boxtimes (\id_{\Gm}-e_1)$ is a
boundary, i.e., there exists $b\in \ZF_m((\Delta^{r+1}\times X)
\times \Gm,Y\times \Gm))$ such that $\partial_i(b)=0$ for
$i=0,1,\dots,r$ and $\partial_{r+1}(b)=a\boxtimes (\id_{\Gm}-e_1)$.
By Lemma \ref{l:colimit} there exists an $N$ such that
$b\in \ZF^{(N)}_m(\Delta^{r+1}\times X \times \Gm,Y\times \Gm )$.
Since
$\ZF^{(N)}_m(-\times \Gm,Y\times \Gm )$
is a presheaf, then
$\partial_i(b)\in \ZF^{(N)}_m(\Delta^{r}\times X \times \Gm,Y\times
\Gm )$. Since $\rho_N$ is a presheaf morphism $\ZF^{(N)}_m(-\times X
\times \Gm,Y\times \Gm )\to \ZF_{m+1}(-\times X,Y)$, one has
$\partial_i(\rho_N(b)) = \rho_N(\partial_i(b))$. Thus,
\begin{multline*}
\partial_i(\rho_N(b)) = \rho_N(\partial_i(b)) = 0 \ \ \text{for} \ \ 0\leq i\leq r, \\
\partial_{r+1}(\rho_N(b)) = \rho_N(\partial_{n+1}(b)) = \rho_N(a\boxtimes
(\id_{\Gm}-e_1)).
\end{multline*}
We see that the homology class of the Moore cycle $\rho_N(a\boxtimes
(\id_{\Gm}-e_1))$ vanishes. By Lemma~\ref{l:Moore_1} the homology
class of the Moore cycle $a\boxtimes \epsilon$ vanishes in
$H_r(\ZF_{m+1}(\Delta^{\bullet}\times X,Y))$. Thus the homology
class of the Moore cycle $a\boxtimes \epsilon\boxtimes \epsilon$
vanishes in
$H_r(\ZF_{m+2}(\Delta^{\bullet}\times X,Y))$. By
Corollary~\ref{box_epsilon_2_and_Sigma_2} the homology class of
$a\boxtimes \Sigma^2$ vanishes in
$H_r(\ZF_{m+2}(\Delta^{\bullet}\times X,Y))$,
too.
\end{proof}

\begin{lemma}\label{l:Moore_2}
Let $X$ and $Y$ be $k$-smooth schemes and $m,r\geq 0$ be integers.
Let $n$ be the integer from Lemma \ref{l:switch2}. Then for any
Moore cycle $a \in \ZF_m((\Delta^r \times X)\wedge (\Gm,1),Y\wedge
(\Gm,1))$ there exists an integer $N$ such that the element
$\rho_N(a)$ is defined and the homology class of the Moore cycle
$$\Sigma^{2n}(\rho_N(a))\boxtimes (id-e_1) \in \ZF_{m+2n+1}((\Delta^r \times X)\wedge (\Gm,1),Y\wedge (\Gm,1))$$
coincides with the homology class of the Moore cycle
$\Sigma^{2n}(a\boxtimes \epsilon)$.
\end{lemma}

\begin{proof}
Set $a'=a\boxtimes (id-e_1)$. Let $H$ be the $\A^1$-homotopy from
Lemma \ref{l:switch2}. Consider the element
$$H(a') \in \ZF_{m+2n}((\Delta^r \times X)\times \Gm\times \Gm,Y\times \Gm\times \Gm).$$
By Lemma \ref{l:colimit} there is an integer $N$ such that
$$a \in \ZF^{(N)}_m((\Delta^r \times X)\times \Gm,Y\times \Gm)$$
and $$H(a') \in \ZF^{(N)}_{m+2n}((\Delta^r \times
X)\times\Gm\times\Gm \times \A^1,Y\times \Gm\times \Gm).$$ Since
$a'$ is a Moore cycle and $H$ is a presheaf morphism, the element
$H(a')$ is a Moore cycle in $\ZF_{m+2n}((\Delta^{\bullet} \times
X)\times \Gm\times \Gm \times \A^1,Y\times \Gm\times \Gm)$. Since
$$\ZF^{(N)}_{m+2n}((-\times X)\times \Gm\times \Gm\times \A^1,Y\times \Gm\times \Gm)$$
is a subpresheaf of $\ZF_m((-\times X)\times \Gm\times \Gm\times
\A^1,Y\times \Gm\times \Gm)$, it follows that $H(a')$ is a Moore
cycle in $\ZF^{(N)}_{m+2n}((\Delta^{\bullet} \times X)\times
\Gm\times \Gm\times \A^1,Y\times \Gm\times \Gm)$.

Applying the presheaf morphism
$$\rho_N: \ZF^{(N)}_{m+2n}((- \times X)\times \Gm\times \Gm\times \A^1,Y\times \Gm\times \Gm) \to
\ZF_{m+2n+1}((- \times X)\times \Gm\times \A^1,Y\times \Gm)$$ to the
Moore cycle $H(a')$, we get a Moore cycle
$$\rho_N(H(a'))\in \ZF_{m+2n+1}((\Delta^r \times X)\times \Gm\times \A^1,Y\times \Gm).$$
Hence $i^*_0(\rho_N(H(a')))\in \ZF_{m+2n+1}((\Delta^r \times
X)\times \Gm,Y\times \Gm)$ and $i^*_1(\rho_N(H(a')))\in
\ZF_{m+2n+1}((\Delta^r \times X)\times \Gm,Y\times \Gm)$ are Moore
cycles, too. Furthermore,
$$i^*_0(\rho_N(H(a')))=\rho_N(i^*_0(H(a')))=\rho_N(\Sigma^{2n}(a'))=\Sigma^{2n}(\rho_N(a'))$$
and
$$i^*_1(\rho_N(H(a')))=\rho_N(i^*_1(H(a')))=\rho_N(\Sigma^{2n}[(\id_Y\times \tau)\circ a' \circ (\id_X\times \tau)]).$$
The two morphisms
$$
i^*_0,i^*_1: \ZF_{n+2m+1}((\Delta^\bullet \times X)\times \Gm\times
\A^1,Y\times \Gm))\rightrightarrows \ZF_{n+2m+1}((\Delta^\bullet
\times X)\times \Gm,Y\times \Gm))
$$
of simplicial abelian groups induce the same morphisms on homology.
The element $\rho_N((H(a'))$ is a Moore cycle. Thus the homological
classes of the Moore cycles $i^*_0(\rho_N(H(a')))$ and
$i^*_1(\rho_N(H(a')))$ coincide in $H_r(\ZF_{n+2m+1}((\Delta^\bullet
\times X)\times \Gm,Y\times \Gm)))$.

By Lemma \ref{rho_and_suspension} one has
$\rho_N(\Sigma^{2n}(a'))=\Sigma^{2n}(\rho_N(a'))$. Thus the first
homological class is the class of $\Sigma^{2n}(\rho_N(a'))=
\Sigma^{2n}(\rho_N(a\boxtimes (id-e_1)))$. By Lemma \ref{l:Moore_1}
the latter homological class coincides with the class of the element
$\Sigma^{2n}(a\boxtimes \epsilon)$.

The element $i^*_1(\rho_N(H(a')))$ coincides with
$\rho_N(\Sigma^{2n}[(\id_Y\times \tau)\circ (a\boxtimes (id-e_1))
\circ (\id_X\times \tau)])$. By Lemma \ref{rho_and_suspension} the
latter element coincides with
$$\Sigma^{2n}(\rho_N[(\id_Y\times \tau)\circ (a\boxtimes (id-e_1)) \circ (\id_X\times \tau)])=\Sigma^{2n}[\rho_N(a)\boxtimes (id-e_1)].$$
Hence the homological classes $\Sigma^{2n}(a\boxtimes \epsilon)$ and
$[\Sigma^{2n}[\rho_N(a)\boxtimes (id-e_1)]]$ coincide in
$H_r(\ZF_{n+2m+1}((\Delta^\bullet \times X)\times \Gm,Y\times
\Gm))$. Finally, the complex $\ZF_{n+2m+1}((\Delta^\bullet \times
X)\wedge (\Gm,1),Y\wedge (\Gm,1))$ is a direct summand in
$\ZF_{n+2m+1}((\Delta^\bullet \times X)\times \Gm,Y\times \Gm)$ and
the elements $\Sigma^{2n}(a\boxtimes \epsilon),
\Sigma^{2n}(\rho_N(a)\boxtimes (id-e_1))$ are in
$\ZF_{n+2m+1}((\Delta^r \times X)\wedge (\Gm,1),Y\wedge (\Gm,1))$.
Hence the homological classes $[\Sigma^{2n}[\rho_N(a)\boxtimes
(id-e_1)]]$ and $[\Sigma^{2n}(a\boxtimes \epsilon)]$ coincide in
$H_r(\ZF_{n+2m+1}((\Delta^\bullet \times X)\wedge (\Gm,1),Y\wedge
(\Gm,1)))$.
\end{proof}

\begin{lemma}\label{l:surjectivity}
Let $X$ and $Y$ be $k$-smooth schemes and $m,r\geq 0$ be integers.
Let $n$ be the integer from Lemma \ref{l:switch2}. Then
\begin{multline*}
\text{Im}[(-\boxtimes \Sigma^{2n+2}): H_r(\ZF_m((\Delta^\bullet \times X)\wedge (\Gm,1),Y\wedge (\Gm,1))) \to \\
\to H_r(\ZF_{m+2n+2}((\Delta^\bullet \times X)\wedge (\Gm,1),Y\wedge (\Gm,1)))]\subseteq \\
\text{Im}[(-\boxtimes (\id_{\Gm}-e_1)):
H_r(\ZF_{m+2n+2}(\Delta^\bullet \times X,Y)) \to
H_r(\ZF_{m+2n+2}((\Delta^\bullet \times X)\wedge (\Gm,1),Y\wedge
(\Gm,1)))].
\end{multline*}
\end{lemma}

\begin{proof}
Take a Moore cycle $a' \in \ZF_m((\Delta^r \times X)\wedge
(\Gm,1),Y\wedge (\Gm,1))$. Then the element $a:=a'\boxtimes
\epsilon$ is a Moore cycle in $\ZF_{m+1}((\Delta^r \times X)\wedge
(\Gm,1),Y\wedge (\Gm,1))$. By Lemma \ref{l:Moore_2} the homology
classes of $\Sigma^{2n}(a\boxtimes \epsilon)$ and
$\Sigma^{2n}(\rho_N(a))\boxtimes (id-e_1)$ coincide in
$$H_r(\ZF_{m+2+2n}((\Delta^{\bullet} \times X)\wedge (\Gm,1),Y\wedge (\Gm,1))).$$
By Corollary \ref{box_epsilon_2_and_Sigma_2} the homology classes of
$\Sigma^{2n}(a\boxtimes \epsilon)=\Sigma^{2n}(a'\boxtimes
\epsilon\boxtimes \epsilon)$ and $\Sigma^{2n+2}(a')$ coincide. Hence
the homology classes of $\Sigma^{2n+2}(a')$ and
$\Sigma^{2n}(\rho_N(a'\boxtimes \epsilon))\boxtimes (id-e_1)$
coincide in $H_r(\ZF_{m+2+2n}((\Delta^{\bullet} \times X)\wedge
(\Gm,1),Y\wedge (\Gm,1)))$.
\end{proof}

We are now in a position to prove Theorem~C.

\begin{theoremc}
Let $X$ and $Y$ be $k$-smooth schemes. Then
\[
-\boxtimes (\id_\Gm - e_1) \colon \ZF(\Delta^\bullet \times X,Y) \to
\ZF((\Delta^\bullet \times X)\wedge (\Gm,1),Y\wedge (\Gm,1))
\]
is a quasi-isomorphism of complexes of abelian groups.
\end{theoremc}

\begin{proof}
The theorem follows from Lemmas \ref{l:injectivity} and
\ref{l:surjectivity}.

In more detail, first prove that the morphism $-\boxtimes (\id_\Gm - e_1)$ induces an epimophism on 
homology groups. 
For this take an integer $r\geq 0$ and an element 
$a \in H_r(\ZF((\Delta^\bullet \times X)\wedge (\Gm,1),Y\wedge (\Gm,1)))$. 
We will find an element 
$b\in H_r(\ZF(\Delta^\bullet \times X,Y))$
such that 
$b\boxtimes (\id_\Gm - e_1)=a$. Note that
there exist an integer $m\geq 0$ and an element 
$a_m\in H_r(\ZF_m((\Delta^\bullet \times X)\wedge (\Gm,1),Y\wedge (\Gm,1)))$
which is a lift of the element $a$.
Let $n$ be the integer from Lemma \ref{l:switch2}. 
By Lemma \ref{l:surjectivity}
there exists an element 
$b_{m+2n+2}\in H_r(\ZF_{m+2n+2}(\Delta^\bullet \times X,Y))$
such that 
$$b_{m+2n+2}\boxtimes (\id_{\Gm}-e_1)=a_m\boxtimes \Sigma^{2n+2}
\in H_r(\ZF_{m+2n+2}((\Delta^\bullet \times X)\wedge (\Gm,1),Y\wedge (\Gm,1))).$$
Let $b$ be the image of $b_{m+2n+2}$ in $H_r(\ZF(\Delta^\bullet \times X,Y))$. Clearly, 
$b\boxtimes (\id_{\Gm}-e_1)=a$ in $H_r(\ZF((\Delta^\bullet \times X)\wedge (\Gm,1),Y\wedge (\Gm,1)))$.
Thus the morphism $-\boxtimes (\id_\Gm - e_1)$ induces an epimophism on homology. 
The fact that the morphism $-\boxtimes (\id_\Gm - e_1)$ induces a monomophism on homology
is proved in a similar fashion. Theorem C is proved.
\end{proof}


\appendix\section{}

The main goal of this section is to prove Theorem D.

Let $(\cc V,\otimes)$ be a closed symmetric monoidal category and
$\cc C$ is a bicomplete category which is tensored and cotensored
over $\cc V$. Then for every $V\in\cc V$ and $C\in\cc C$ there are
defined objects $V\otimes C$, $C\otimes V$, $\underline{\Hom}(V,C)$
of $\cc C$. They are all functorial in $V$ and $C$. Moreover, for
every morphism $r:V\to V'$ in $\cc V$ the square in $\cc C$
   \begin{equation}\label{actionV}
    \xymatrix{C\ar[rr]^{-\otimes V}\ar[d]_{-\otimes V'}&&\underline{\Hom}(V,C\otimes V)\ar[d]^{r_*}\\
    \underline{\Hom}(V',C\otimes V')\ar[rr]_{r^*}&&\underline{\Hom}(V,C\otimes V')}
   \end{equation}
is commutative.

{As an important example}, $\cc V$ is
the category of simplicial objects $sPre(\ZF_0(k))$ in the category
$Pre(\ZF_0(k))$ and $\cc C$
is the category
$sPre_{Ab}(\ZF_*(k))$
of simplicial objects in $Pre_{Ab}(\ZF_*(k))$.
The General Framework of p.~\pageref{gfr} is then immediately extended to this couple $(\cc V,\cc C)$.
Recall that the functor $\ZF_*(k)\times \ZF_0(k) \xrightarrow{\boxtimes} \ZF_*(k)$
takes $(X,Y)$ to $X\times Y$. As usual, the Yoneda embedding identifies the 
category simplicial objects in $\ZF_0(k)$  
with a full subcategory of $sPre_{Ab}(\ZF_0(k))$.


The following lemma is obvious.

\begin{lemma}\label{l:3_of_4}
Suppose in the diagram \eqref{actionV} the morphisms $r_*$, $r^*$ and $-\otimes V'$
are sectionwise weak equivalences, then the morphism $-\otimes V$
is a sectionwise weak equivalence.
\end{lemma}

As it is shown in~\cite[Section~5]{GP1}, the category of framed correspondences
of level zero $\Fr_0(k)$ has an action by finite pointed sets $S\otimes
K:=\bigsqcup_{K\setminus *}S$ with $S\in Sm/k$ and $K$ a finite pointed
set. The cone of $S$ is the simplicial object $S\otimes I$ in
$\Fr_0(k)$, where $(I,1)$ is the pointed simplicial set $\Delta[1]$
with basepoint 1. There is a natural morphism $i_0:S\to S\otimes I$
in $\Delta^{\textrm{\op}}\Fr_0(k)$.
Let $\pt\xrightarrow{e_1} \Gm$ be the point $1$ in $\Gm(k)$. Then
$\bb G_m^{\wedge 1}$ is the simplicial object in $\Fr_0(k)$
which is obtained by taking the pushout of the diagram
   $\Gm\xleftarrow{e_1} \pt\bl{i_0}\hookrightarrow \pt\otimes I$
in $\Delta^{\textrm{\op}}\Fr_0(k)$.

Let $L: \Fr_0\to \ZF_0$ be the canonical functor which is the identity on objects and which takes
a morphism $\phi\in \Fr_0(Y,X)$ to the class $1\cdot \phi\in \ZF_0(Y,X)$.
If we apply the functor $L$ to $\bb G_m^{\wedge 1}$, we get an object in
$sPre_{Ab}(\ZF_0(k))$ denoted by $\ZF_0(\bb G_m^{\wedge 1})$.

Put $\ZF_0(\bb G_m,1)=\ZF_0(\bb G_m)/Im(e_{1,*})=\text{Ker}(e^*_1)$.
There is a unique morphism $r: \ZF_0(\bb G_m^{\wedge 1})\to \ZF_0(\bb G_m,1)$
which restricts to the quotient map $q: \ZF_0(\bb G_m)\to \ZF_0(\bb G_m)/Im(e_{1,*})$
on $\ZF_0(\bb G_m)$ and which restricts to the zero map on $\pt\otimes I$.

The following lemma is straightforward and left to the reader.

\begin{lemma}
$\ZF(X)\boxtimes \ZF_0(\bb G_m,1)=\ZF(X\wedge (\bb G_m,1))$,
$\ZF(X)\boxtimes \ZF_0(\bb G_m^{\wedge 1})=\ZF(X\times \bb G_m^{\wedge 1})$.
\end{lemma}

\begin{lemma}\label{l:r_and_r}
The morphisms
$$r_*: \underline{\Hom}(\ZF_0(\bb G_m^{\wedge 1}),C_*(\ZF(X)\boxtimes \ZF_0(\bb G_m^{\wedge 1})))\to
\underline{\Hom}(\ZF_0(\bb G_m^{\wedge 1}),C_*(\ZF(X)\boxtimes \ZF_0(\bb G_m,1))))$$
$$r^*: \underline{\Hom}(\ZF_0(\bb G_m,1),C_*(\ZF(X)\boxtimes \ZF_0(\bb G_m,1)))\to
\underline{\Hom}(\ZF_0(\bb G_m^{\wedge 1}),C_*(\ZF(X)\boxtimes \ZF_0(\bb G_m,1))))$$
are sectionwise weak equivalences.
\end{lemma}

\begin{proof}
It is easy to see that the morphisms
$$r: \ZF_0(\bb G_m^{\wedge 1})\to \ZF_0(\bb G_m,1),$$
$$id\boxtimes r: \ZF(X)\boxtimes \ZF_0(\bb G_m^{\wedge 1})\to \ZF(X)\boxtimes \ZF_0(\bb G_m,1),$$
$$id\boxtimes r: C_*(\ZF(X)\boxtimes \ZF_0(\bb G_m^{\wedge 1}))\to C_*(\ZF(X)\boxtimes \ZF_0(\bb G_m,1))$$
are sectionwise weak equivalences. The lemma now follows.
\end{proof}

Theorem C, Lemma \ref{l:3_of_4} and Lemma \ref{l:r_and_r} imply the following

\begin{corollary}\label{cor:_C1}
The morphism
$$- \boxtimes \bb G_m^{\wedge 1}: C_*\ZF(X)\to \underline{\Hom}(\bb G_m^{\wedge 1},C_*\ZF(X\times \bb G_m^{\wedge 1}))$$
is a sectionwise weak equivalence in $sPre_{Ab}(\ZF_*(k))$.
\end{corollary}

We are now in a position to prove the following

\begin{theoremd}\label{thm:D}
The morphism $c_0: LM_{fr}(X)\to \underline{\Hom}(\bb G,LM_{fr}(X\times\bb G_m^{\wedge 1}))$
is a sectionwise stable weak equivalence of presheaves of $S^1$-spectra.
\end{theoremd}

\begin{proof}
First, the adjunction unit
$\bb G\xrightarrow{\textrm{adj}} (\bb G_m^{\wedge 1}|_{Sm/k})$ in $sPre_\bullet(Sm/k)$
induces an isomorphism
$
\underline{\Hom}(\bb G_m^{\wedge 1},LM_{fr}(X\times\bb G_m^{\wedge 1}))\xrightarrow{\textrm{adj}^*}
\underline{\Hom}(\bb G,LM_{fr}(X\times\bb G_m^{\wedge 1}))
$
of $S^1$-spectra.
Second, the morphism $\textrm{adj}^*\circ (- \boxtimes \bb G_m^{\wedge 1})$
coincides with
the morphism
$$c_0: LM_{fr}(X)\to \underline{\Hom}(\bb G,LM_{fr}(X\times\bb G_m^{\wedge 1})).$$
is the morphism~\eqref{strelki_c}. The theorem now follows from Corollary~\ref{cor:_C1}.
\end{proof}

\section{Some facts on henzelization}\label{henzelization}

We refer the reader to~\cite{Gab} or~\cite{FP} for the definition of henzelization
of an affine scheme along a closed subscheme.

Let $X,X_1$ be $k$-smooth affine varieties, $Z\subset X$, $Z_1\subset X_1$ be closed subsets.
Let $f: X_1\to X$ be a $k$-morphism such that $Z_1\subset f^{-1}(Z)$.
For an \'{e}tale neighborhood $(W,\pi: W\to X,s:Z\to W)$ of $Z$ in $X$
set $W_1=X_1\times_X W$. Let $\pi_1: W_1\to X_1$ be the projection and let
$s_1=(i_1,f|_{Z_1)}: Z_1 \to W_1$, where $i_1: Z_1 \hookrightarrow X_1$ be the inclusion.
Then $(W_1,\pi_1,s_1)$ is an \'{e}tale neighborhood of $Z_1$ in $X_1$.
Denote by $f_W: W_1\to W$ the projection. Then one has a morphism
$\lim (f_W): \lim_{(W,\pi,s)} W_1 \to \lim_{(W,\pi,s)} W=X^h_Z$. Set,
   \begin{equation}\label{f_h}
    f^h=\lim (f_W)\circ can_f: (X_1)^h_{Z_1} \to X^h_Z,
   \end{equation}
where $can_f: (X_1)^h_{Z_1}\to \lim_{(W,\pi,s)} W_1$ is the canonical morphism.
Clearly, $\rho\circ f^h=f\circ \rho_1$, where  $\rho: X^h_Z\to X$ and
$\rho_1: (X_1)^h_{Z_1}\to X_1$ are the canonical morphisms.

The following properties of the morphism $f^h$ are straightforward:

\begin{enumerate}
\item For any affine $k$-smooth variety $X$ one has $\id^h_X=\id_{X^h_Z}$. If $p: X\to \pt$
is the structure map, then for any closed $Z$ in $X$ the morphism
$p^h: X^h_Z \to (\pt)^h_{\pt}=\pt$
is the structure morphism.

\item Given a $k$-morphism $f_1: X_2\to X_1$ of affine $k$-smooth varieties and
a closed subset $Z_2\subset X_2$ with $Z_2\subset f^{-1}_1(Z_1)$
one has
$(f\circ f_1)^h=f^h\circ f^h_1$.


\item If $i: Z\hookrightarrow X$ is the closed inclusion, $Z_1=Z$, then
$Z^h_Z=Z$ and $i^h: Z=Z^h_Z\to X^h_Z$ coincides with the canonical closed inclusion
$s: Z\to X^h_Z$.
\end{enumerate}

The last two properties imply the following property.
Let $X$ be an affine $k$-smooth variety and $x\in X$ be a $k$-rational point.
Suppose $s: \pt \to X^h_x$ is the closed point of $X^h_x$ and $i_x: \pt \to X$ is the point $x$.
Let $p: X\to \pt$ be the structure map. Then one has equalities
   $$(i_x\circ p)^h=i^h_x \circ p^h=s\circ p^h: X^h_x \to X^h_x.$$

These observations imply the following

\begin{lemma}\label{f_h_s}
Let $X$ be an affine $k$-smooth variety and $x\in X$ be its $k$-rational point.
Let $f_s: \A^1\times X \to X$ be a morphism such that $f_1: X\to X$ is the identity,
$f_0: X\to X$ coincides with the morphism $X\xrightarrow{p} \pt\xrightarrow{i_x} X$
and $f_s(\A^1\times \{x\})=\{x\}$. Then the morphism
$f^h_s: (\A^1\times X)^h_{\A^1\times x} \to X^h_x$
defined by~\eqref{f_h} has the following properties:
\begin{enumerate}
\item $(f^h_s)|_{(1\times X)^h_{(1,x)}}: X^h_x\to X^h_x$ is the identity;
\item $(f^h_s)|_{(0\times X)^h_{(0,x)}}: X^h_x\to X^h_x$ coincides with the morphism $X^h_x\xrightarrow{p^h} \pt \xrightarrow{s_x} X^h_x$,
where $s_x: \pt \hookrightarrow X^h_x$ is the closed point of $X^h_x$.
\end{enumerate}
\end{lemma}

\begin{proof}
The first assertion follows from the equalities
   $$\id_{X^h_Z}=\id^h_X=(f_1)^h=(f_s\circ i_1)^h=f^h_s \circ i^h_1=(f^h_s)|_{(1\times X)^h_{(1,x)}}.$$
Let $s_x: \pt \hookrightarrow X^h_x$ be the closed point of $X^h_x$.
As mentioned above, $s_x=i^h_x$, where $i_x: \pt \xrightarrow{} X$ is the closed point $x$ of $X$.
The equalities
   $$s_x\circ p^h=i^h_x \circ p^h=(i_x\circ p)^h=f^h_0=(f_s\circ i_0)^h=f^h_s \circ i^h_0=(f^h_s)|_{(0\times X)^h_{(0,x)}}$$
imply the second assertion.
\end{proof}

If we take $X=\A^m$, a $k$-rational point $x\in \A^m$ and the morphism
$f_s: \A^1\times \A^m \to \A^m$ mapping $(s,y)$ to $s\cdot (y-x)+x$,
then $f_s: \A^1\times X \to X$ satisfies
the hypotheses of Lemma~\ref{f_h_s}.
Thus Lemma~\ref{f_h_s} implies the following statement, which is in fact Lemma~\ref{Retraction}.

\begin{corollary}
The morphism $H_{s}:=f^h_s: U'_{s}\to U'$ has the following properties:
\begin{itemize}
\item[(a)] $H_1:=(f^h_s)|_{(1\times X)^h_{(1,x)}}: U' \to U'$ is the identity morphism;
\item[(b)] $H_0:=(f^h_s)|_{(0\times X)^h_{(0,x)}}: U' \to U'$ coincides with the composite morphism
$U'\xrightarrow{p^h} pt \xrightarrow{s_x} U'$, where
$p^h: U' \to pt=\spec(k)$ is the structure morphism and
$s_x: \pt \hookrightarrow X^h_x$ is the closed point of $X^h_x$.
\end{itemize}
\end{corollary}


\section*{Acknowledgements}
This paper was partly written during the visit of the authors in
summer 2014 to the University of Duisburg-Essen (Marc Levine's
Arbeitsgruppe). They would like to thank the University for its kind
hospitality and support. The first author and the third author were 
supported by a grant from the Government of the Russian Federation (agreement 
075-15-2019-1620). The first author was also supported by the Young Russian Mathematics award.


\begin{thebibliography}{XXXX}

\bibitem[DP]{DP}
A.~Druzhinin, I.~Panin, \emph{Surjectivity of the etale excision map for homotopy invariant framed presheaves},
preprint arXiv:1808.07765.

\bibitem[FP]{FP}
R.~Fedorov, I.~Panin, \emph{A proof of
the Grothendieck--Serre conjecture on principal bundles over regular local rings containing infinite fields},
Publ. Math. IHES 122(1) (2015), 169-193.

\bibitem[Gab]{Gab}
O.~Gabber, \emph{Affine analog of the proper base change theorem}, Israel J. Math., 87 (1994), 325-335.

\bibitem[GNP]{GNP}
G.~Garkusha, A.~Neshitov, I.~Panin, \emph{Framed motives of relative motivic spheres}, Trans. Amer.
Math. Soc., to appear. DOI https://doi.org/10.1090/tran/8386.

\bibitem[GP1]{GP1}
G.~Garkusha, I.~Panin, \emph{Framed motives of algebraic varieties (after V. Voevodsky)}, 
J. Amer. Math. Soc. 34(1) (2021), 261-313.

\bibitem[GP2]{GP2} G.~Garkusha, I.~Panin, \emph{Homotopy invariant presheaves with framed transfers},
Cambridge J. Math. 8(1) (2020), 1-94.

\bibitem[GP3]{GP2a} G.~Garkusha, I.~Panin, \emph{The triangulated categories of framed bispectra and framed motives}, preprint
arXiv:1809.08006.

\bibitem[GP4]{GP3} G. Garkusha, I. Panin, \emph{On the motivic spectral sequence},
J. Inst. Math. Jussieu 17(1) (2018), 137-170.

\bibitem[EGA4]{LNM146} A. Grothendieck, \emph{\'El\'ements de g\'eom\'etrie alg\'ebrique. IV, \'Etude
locale des schemas et des morphismes de schemas}, Quatri\`eme
partie, Publ. Math. IHES 32 (1967), 5-361.

\bibitem[GrD]{GrD} A. Grothendieck, J. Dieudonn\'{e}, \'El\'{e}ments de G\'{e}om\'{e}trie Alg\'{e}brique
IV. \'Etude locale des sch\'{e}mas et des morphismes de sch\'{e}mas (Troisi\`{e}me
Partie), Publ. Math. IH\'ES 28 (1966), 5-255.

\bibitem[Jar]{Jar} {J. F. Jardine}, \emph{Motivic symmetric spectra}, Doc. Math. 5 (2000), 445-552.

\bibitem[Hir]{Hir}{Ph. S. Hirschhorn}, \emph{Model categories and their
         localizations}, Mathematical Surveys and Monographs 99, 2003.

\bibitem[MV]{MV} F. Morel, V. Voevodsky, $\bb A^1$-homotopy theory of schemes,
         Publ. Math. IHES 90 (1999), 45-143.

\bibitem[OP]{OP} M.Ojanguren, I.Panin, Rationally trivial hermitian spaces are locally trivial,
Math. Z. 237(1) (2001), 181-198.

\bibitem[Sch]{Sch} S. Schwede, \emph{An untitled book project about symmetric spectra},
         available at www.math.uni-bonn.de/people/schwede/SymSpec-v3.pdf (version April 2012).

\bibitem[Seg]{Seg} {G. Segal}, \emph{Categories and cohomology theories}, Topology 13 (1974), 293-312.

\bibitem[S]{S}
A.~Suslin, \emph{On the Grayson spectral sequence}, Proc. Steklov
Inst. Math. 241 (2003), 202-237.

\bibitem[V1]{V2} V.~Voevodsky, \emph{Notes on framed correspondences}, unpublished, 2001. Also available at
https://www.math.ias.edu/vladimir/publications

\bibitem[V2]{V1}
V.~Voevodsky, \emph{Cancellation theorem}, Doc. Math. Extra Volume in honor of A. Suslin (2010), 671-685.

\end{thebibliography}
\end{document}